\documentclass[onefignum,onetabnum]{siamart220329}

\usepackage{amsfonts}
\usepackage{graphicx}
\usepackage{epstopdf}
\usepackage{algorithmic}
\ifpdf
\DeclareGraphicsExtensions{.eps,.pdf,.png,.jpg}
\else
\DeclareGraphicsExtensions{.eps}
\fi

\newsiamremark{remark}{Remark}
\newsiamremark{hypothesis}{Hypothesis}
\crefname{hypothesis}{Hypothesis}{Hypotheses}
\newsiamthm{claim}{Claim}

\headers{A Fully Parallelized and Budgeted MLMC Method}{N. Baumgarten, S. Krumscheid, and C. Wieners}

\title{
    A Fully Parallelized and Budgeted Multi-level Monte Carlo Method
    and the Application to Acoustic Waves\thanks{Submitted to the editors \today}
    \funding{This work was partly funded by the German Research Foundation (DFG) through an association with
      the CRC 1173 on Wave Phenomena
      (Project-ID 258734477). The numerical experiments were executed on the
    Hochleistungsrechner Karlsruhe (HoReKa) funded by the
    Ministry of Science, Research and the Arts Baden-Württemberg and by the Federal Ministry of Education
    and Research.}
}

\author{
    Niklas Baumgarten\thanks{Institut für Angewandete und Numerische Mathematik, KIT, \email{niklas.baumgarten@kit.edu}.}
    \and Sebastian Krumscheid\thanks{Steinbuch Computing Center (SCC), KIT, \email{sebastian.krumscheid@kit.edu}.}
    \and Christian Wieners\thanks{Institut für Angewandete und Numerische Mathematik, KIT, \email{christian.wieners@kit.edu}.}
}

\usepackage{amsopn}

\RequirePackage{mathtools}
\RequirePackage{adjustbox}
\RequirePackage{stmaryrd}
\RequirePackage{graphicx}
\RequirePackage{enumitem}
\RequirePackage{ulem}
\RequirePackage{mdframed}
\RequirePackage{amssymb}
\RequirePackage{pbox}
\input{preamble.sty}

\ifpdf
\hypersetup{
    pdftitle={A Fully Parallelized and Budgeted Multi-level Monte Carlo Method
    and the Application to Acoustic Waves},
    pdfauthor={N. Baumgarten, S. Krumscheid, and C. Wieners}
}
\fi

\begin{document}

    \maketitle

    \begin{abstract}
    We present a novel variant of the multi-level Monte
    Carlo method that effectively utilizes a reserved computational
    budget on a high-performance computing system to minimize the mean
    squared error. Our approach combines concepts of the continuation
    multi-level Monte Carlo method with dynamic programming techniques
    following Bellman's optimality principle,
    and a new parallelization strategy based on a single distributed data
    structure. Additionally, we establish a theoretical bound on the error
    reduction on a parallel computing cluster and provide empirical
    evidence that the proposed method adheres to this
    bound. We implement, test, and benchmark
    the approach on computationally demanding problems, focusing on
    its application to acoustic wave propagation in high-dimensional
    random media.
\end{abstract}



    \section{Introduction}\label{sec:introduction}

Being certain about an outcome of any physical, technical or
economical process comes with a cost. This cost manifests in several
ways and often includes conducting extensive research, data
collection, analysis, and employing sophisticated modeling or
simulation techniques. Furthermore, achieving certainty may involve
dealing with extremely complex systems in high dimensions and
intricate mathematical models which may require the usage of highly
advanced and resource-intensive computational technologies.

In this paper, we investigate the computational aspects associated
with the cost of certainty, specifically the cost involved in
quantifying uncertainty. To address this problem, we propose an
integrated framework that combines multi-level Monte Carlo (MLMC),
finite element (FE), and dynamic programming (DP) methods. In
particular, we introduce a fully parallelized and budgeted variant of
the multi-level Monte Carlo method, named as Budgeted MLMC (BMLMC)
method. Additionally, we present novel parallelization concepts
to handle substantial computational loads.

In recent years, the combination of MLMC and FE methods has been
successfully applied to partial differential equations (PDEs)
involving random parameters. Notable instances include elliptic
PDEs~\cite{abdulle2013multilevel,
    barth2011multi, charrier2012strong, charrier2013finite,
    cliffe2011multilevel, collier2015continuation,
    teckentrup2013further} as well as hyperbolic
PDEs~\cite{ballesio2019multilevel,
    grote2021uncertainty, mishra2012sparse, mishra2012multi,
    mishra2012multilevel, mishra2013multi, mishra2016multi}. The MLMC-FE
method provides significant computational efficiency, enabling the
achievement of desired accuracy levels while considerably reducing
computational requirements compared to a single-level Monte Carlo (MC)
method. Despite the existence of alternative uncertainty
quantification (UQ) methods, such as stochastic collocation
(SC)~\cite{babuvska2007stochastic, babuvska2010stochastic, nobile2008sparse, nobile2009analysis},
quasi-Monte Carlo (QMC) techniques~\cite{caflisch1998monte, graham2011quasi, sloan1998quasi}, as
well as various multi-level and multi-index variants described
in~\cite{haji2016multiindexsc, haji2016multi, kuo2012quasi, teckentrup2015multilevel},
the MLMC method is a popular choice due
to its non-intrusive implementation and the moderate assumptions it
imposes on the problem and the employed discretization.

The problems that can be addressed using the method proposed in this paper fall
under the domain of forward UQ methods. These problems typically involve
solving linear or nonlinear partial differential equations (PDEs) with
random input data. In our study, we specifically concentrate on the
application of the method to acoustic wave equations in random
media. This particular model is of high significance in applications such as
seismic imaging and geophysics. As a hyperbolic PDE, it presents
unique challenges, including issues related to low regularity,
stability conditions, and the handling of high-dimensional random
input data.

To tackle this challenge, we developed a high-performance computing (HPC)
approach and drew inspiration from previous works such
as~\cite{mishra2016multi, vsukys2013multi}, which address similar
problems using MLMC methods. Additionally, we refer to studies such
as~\cite{bohlen2020visco, corallo2022space, dorfler2016space,
    hochbruck2015efficient}, which focus on different discretization
techniques for this particular problem. However, the size of the
resulting discrete model, the number of required samples, and the
diverse array of solution approaches prompted us to explore new
strategies for distributing the computational workload and effectively
managing the interplay of all algorithms involved. Consequently, our
research and contributions can be summarized by the following three
points.

\paragraph{Budgeted multi-level Monte Carlo method}

We present a novel budgeted variant of the MLMC method which we call
BMLMC method. In the classical MLMC framework described
in~\cite{barth2011multi, cliffe2011multilevel,
    charrier2013finite, giles2008improved, giles2008multilevel,
    teckentrup2013further}, a desired tolerance for the total root mean
squared error (RMSE), denoted as $\epsilon > 0$, is
selected. Subsequently, a sequence of samples is generated either
based on knowledge of the solution's regularity or by adaptive
methods, such as the Continuation MLMC (CMLMC)
method~\cite{collier2015continuation, giles2015multilevel}. The objective
of the adaptive methods is minimizing the overall computational cost
while achieving the specified tolerance. In contrast, the BMLMC method
is designed to minimize the RMSE and to operate within a given cost
budget denoted as $\rB > 0$. The cost and the budget are specifically
measured in units of CPU seconds. A similar approach is also taken by
the multi-fidelity method~\cite{peherstorfer2018survey}
and the unbiased estimation in~\cite{rhee2015unbiased}. This point of
view is motivated in the context of this work by two key reasons.

In HPC applications, it is customary to
allocate a cost budget to secure a spot in the queue of an HPC
cluster. In our notation, this is represented as $\rB = \abs{\cP}
\cdot \rT_{\mathrm{\rB}}$, where $\abs{\cP}$ corresponds to the total
number of processing units, and $\rT_{\mathrm{\rB}}$ represents the
time budget. Even with a priori knowledge about the solution's
regularity, the convergence rates, the parallelization strategies and
the computational infrastructure, determining the tolerance $\epsilon$
that can be achieved with a budget $\rB$ remains an NP-hard
task. Thus, it is more natural in HPC applications to replace the
predefined tolerance $\epsilon$ with a predefined budget $\rB$ and
allow the method to determine the smallest achievable $\epsilon$ while
fully utilizing the budget. This approach aligns with the practical
requirements and constraints encountered in HPC scenarios.

Furthermore, employing a budgeted algorithm enables us to address the
question of determining the optimal algorithm stack for a given
problem. Specifically, we aim to identify the algorithm stack that
yields the smallest error while utilizing the same computational
resources. We can investigate this question empirically by conducting
experiments where various algorithms, e.g.~different time stepping methods,
are evaluated using an equal budget, and subsequently
comparing the resulting error estimates. This
empirical analysis allows us to make informed decisions regarding the
selection of algorithm combinations that optimize the trade-off
between computational resources and the achieved accuracy.

In order to implement our method, we introduce the reserved
computational budget as an additional constraint, effectively
transforming the problem into a knapsack problem. To solve this
knapsack problem, we decompose it into multiple
subproblems, each of which is solved under an optimality condition.
Specifically, we combine the Continuation MLMC (CMLMC)
method~\cite{collier2015continuation} with dynamic programming
techniques, allowing us to address the problem without any prior
knowledge of its specific characteristics or computational
complexity. The method proves to be highly robust and performant, and
can be widely applied to a general class of PDEs in conjunction with
the proposed parallelization strategy.

\paragraph{Parallelization Strategy}

Parallelization is an omnipresent challenge in various applications,
including manufacturing, logistics, and computer networks. The load
distribution of MLMC methods has already been addressed in prior work, such
as~\cite{badia2021massively, drzisga2017scheduling, vsukys2011static,
    vsukys2013adaptive}. However, our parallelization in this
study offers a distinguishing feature compared to most of the previous
work by being done on a single distributed data structure.
This technical choice allows for a more efficient and adaptable load
distribution.

In our proposed parallelization scheme, we leverage a multi-mesh
parallelism on distributed memory, which can be applied not only to
arbitrary FE spaces but also to other non-intrusive UQ
algorithms beyond MLMC. The key idea is to assemble large algebraic
systems on distributed memory across multiple computing nodes, while
dynamically adapting the system's structure based on the number of
samples that need to be computed. This integrated approach seamlessly
combines FE with UQ methods, resulting in a
highly efficient implementation that fully exploits hardware
capabilities without compromising the non-intrusive nature of UQ
methods. We refer to this parallelization strategy throughout the
paper as multi-sample finite element method (MS-FEM).

By adopting this parallelization strategy, we achieve a highly
adaptive and hardware proximal implementation, providing significant
computational benefits by minimizing the communication overhead and processor idling.
In summary, our parallelization approach offers a novel
perspective on distributing computational load and demonstrates its
effectiveness in enhancing the performance of UQ algorithms, including
the MLMC method.

\paragraph{Numerical Experiments and Software}
Lastly, this paper provides a concise insight into the developed
software, M++~\cite{wieners2023mpp315}, along with a range of automated
numerical experiments. An essential aspect of our implementation is
its flexibility, as the methodology can be seamlessly integrated with
arbitrary FE and other non-intrusive UQ
methods. Consequently, our approach allows for the unified application
of these methods to elliptic, parabolic, and hyperbolic
PDEs~\cite{baumgarten2023fully, baumgarten2021parallel}, although our
focus in this paper is specifically on the application of the BMLMC
method to the acoustic wave equation.

It is worth noting that the empirical investigations concerning the
algorithm stack have been meticulously conducted using a fully
automated approach via a continuous delivery pipeline,
ensuring the reproducibility and potential
improvement of the presented results. In particular, we describe
experiments performed within our framework that encompass the
parallelization, model assessment, and methodology evaluation.

\paragraph{Outline}

The paper is structured as follows.
In~\Cref{sec:multi-level-estimators} we establish the
notation, state the underlying assumptions, and present the classical
MLMC method. Subsequently, we derive the BMLMC method using dynamic
programming techniques. Notably, we also demonstrate in this section
that any parallel and adaptive implementation of an MLMC method inherently
contains an error contribution that cannot be eliminated by increasing the
number of processing units. \Cref{sec:multi-level-estimators} is designed to
be independent of the specific application, making it applicable to
other problems as well.
In~\Cref{sec:multi-sample-finite-element-method}, we discuss
details of the multi-sample finite element method (MS-FEM), which
represents an ideal fit for the requirements of the BMLMC method,
offering a highly adaptive and efficient parallelization strategy. We
provide a explanation of how we define a multi-sample
discontinuous Galerkin (dG) finite element space, which serves as the
basis for the subsequent discussion on the discretization of the
acoustic wave equation under uncertainty
in~\Cref{sec:discretization-of-the-acoustic-wave-equation}.
Specifically, we introduce a semi-discrete dG system solved with implicit
time-stepping methods in the time domain.
\Cref{sec:mlmc-acoustic-experiments} presents numerical
results achieved using the proposed framework.
We conclude with a comprehensive discussion and an outlook on further work
in~\Cref{sec:discussion-and-outlook}.

    \section{A Budgeted Multi-level Monte Carlo Method}
\label{sec:multi-level-estimators}

In this section, we present the budgeted multi-level Monte Carlo (BMLMC)
method along with novel approaches for determining the optimal
load distribution among a given set of processing units.
We begin by establishing the notation and outlining the main assumptions
on the stochastic model in \Cref{subsec:assumptions-and-notation}.
Next, in \Cref{subsec:introduction-to-multi-level-monte-carlo},
we provide an overview of the classical multi-level Monte Carlo
method~\cite{giles2008improved, giles2008multilevel}.
\Cref{subsec:implementation-techniques} focuses on implementation
techniques~\cite{collier2015continuation, giles2015multilevel}
that are utilized in conjunction with dynamic programming (DP)
to develop the BMLMC method,
as described in~\Cref{subsec:budgeted-multi-level-monte-carlo}.
Finally, in~\Cref{subsec:parallelization-techniques},
we explain the resulting algorithm as a distributed state machine,
leveraging adaptive parallelization techniques and show
that the error of any parallel and adaptive implementation
of the MLMC method obeys a bound with respect
to the computing time and the amount of processing units.

\subsection{Assumptions and Notation}
\label{subsec:assumptions-and-notation}

In the following, we consider a bounded polygonal domain $\cD \subset \RR^D$
with spatial dimension $D \in \set{1, 2, 3}$,
and a probability space $(\Omega, \mathcal{F}, \mathbb{P})$.
We are interested in solving PDEs with uncertainly determined input data, i.e.,
models of the form $\cL(\omega, \bx) \bu(\omega, \bx) = \bb(\omega, \bx)$,
where $\omega \in \Omega$ corresponds to a specific outcome in the probability space
and $\bu(\omega, \bx)$ represents the solution of the PDE at spatial location $\bx \in \cD$.
The differential operator
$\mathcal{L}(\omega, \bx)$ depends on the input data
and acts on the solution $\bu(\omega, \bx)$, while $\bb(\omega, \bx)$ is a forcing term on the PDE.
Later in this work, we also consider time dependent PDE models, which however,
relax to the above formulation for a fixed time point. Thereby, the time dependence is
neglected throughout this section.
For each $\omega \in \Omega$, the solution lies in a separable Hilbert space $V$.
Additionally, we consider a bounded functional $\rQ(\omega) \coloneqq \rQ(\bu(\omega, \bx))$
that represents a quantity of interest~(QoI) of the solution.

As a start, we consider the objective of the method as to estimate the expected value of $\rQ$
with a prescribed level of accuracy by approximating the PDE using a FE solution $\bu_\ell(\omega) \in V_{\ell}$
and employing an MC method to approximate the expectation.
In this context, $V_{\ell}$ represents a discrete FE space at level $\ell$ and $\bu_{\ell}(\omega, \bx)$
solves the discrete problem $\cL_{\ell}(\omega, \bx)\bu_{\ell}(\omega, \bx) = \bb_{\ell}(\omega, \bx)$.
A specific study of such a system is given
in~\Cref{sec:discretization-of-the-acoustic-wave-equation}
in form of the acoustic wave equation discretized with non-conforming dG elements.
Furthermore, we denote by $\rQ_\ell(\omega) \coloneqq \rQ_\ell(\bu_\ell(\omega, \bx))$ the quantity of interest (QoI)
defined on $V_\ell$ and give samples the index $m$, i.e.,
$\bu_{\ell}(\omega^{(m)}, \bx)$ represents the FE solution computed
using input data corresponding to $\omega^{(m)}$.
To numerically represent the input, we make the following assumption.

\clearpage

\begin{assumption}
    \label{assumption:fdna}
    We rely on the \textit{finite dimensional noise assumption}
    \text{\textrm{(FDNA)}}: the space of outcomes $\Xi$ of any random field $\bY
    \colon \Omega \times \cD \rightarrow \Xi$ is of finite dimension~$K$,
    i.e., any sample can be represented by a vector $\by^{(m)}
    = (y_1, \dots, y_K)^\top \in \Xi \subset \RR^K$.
\end{assumption}

By the FDNA, we can express samples of the input data by the vector
$\by^{(m)} \in \RR^K$, so that we simply write $\bu_\ell^{(m)} \coloneqq \bu_\ell(\by^{(m)}, \bx)$
and $\rQ_\ell^{(m)} \coloneqq \rQ_\ell(\by^{(m)})$.

A Monte Carlo method for a FE solution estimates
the expected value of a QoI
\begin{equation}
    \EE[\rQ]
    \coloneqq \int_{\Omega} \rQ(\omega) \rd \PP
    \approx \frac{1}M \sum_{m=1}^M \rQ_\ell(\by^{(m)})
    \eqqcolon \widehat{\rQ}_{\ell}
    \,,
    \label{eq:mc-estimator}
\end{equation}
depending on independent and identically distributed (iid) samples
$\by^{(m)} \in \Xi \subset \RR^K$ drawn from the distribution of the input data.
Computing $\rQ_\ell(\by^{(m)})$ requires to solve a PDE with a FE
method and thereby is a costly and inexact evaluation admitting a
discretization error $\err_{\text{disc}}$. Thus, estimating the
expected value comes with an estimator bias induced by the FE method
giving a mean squared error (MSE)
\begin{equation}
    \label{eq:mse-mc}
    \mathrm{err}_{\mathrm{MSE}}
    = \EE \squarelr{\roundlr{\widehat{{\rQ}}_{\ell} - \EE[{\rQ}]}^2}
    = \underbrace{M^{-1}  \VV[{\rQ}_\ell]}_{\text{Estimator variance}}
    + \quad {\underbrace{\EE[{\rQ}_\ell - {\rQ}]}_{\text{Bias}}}^2,
\end{equation}
where $\VV[{\rQ}_\ell]$ is the variance of the random variable $\rQ_\ell$
and the corresponding root mean squared error (RMSE) is given by
$\mathrm{err}_{\mathrm{RMSE}} = \sqrt{\mathrm{err}_{\mathrm{MSE}}}$.
To link the admitted error and the computational cost, we introduce
the following definition of a cost-measure.
\begin{definition}
    [$\epsilon$-time, $\epsilon$-cost and cost-measure]
    For given $\epsilon > 0$, the computing time to achieve a root mean squared error
    $\mathrm{err}_{\mathrm{RMSE}} \leq \epsilon$
    is the \textit{$\epsilon$-time} $\rT_{\epsilon}$.
    The \textit{$\epsilon$-cost} $\rC_{\epsilon}= \abs{\cP} \cdot \rT_\epsilon$
    is the corresponding computational cost on a parallel machine,
    where $\abs{\cP}$ is the total count of involved processing units.
\end{definition}

On a serial machine the $\epsilon$-cost simplifies to
$\rC_{\epsilon,\rs} = \rT_\epsilon$. Since the classical theory
of MLMC methods does not consider parallel machines, the methods
in subsections~\ref{subsec:introduction-to-multi-level-monte-carlo}--\ref{subsec:implementation-techniques}
are formulated for a serial $\epsilon$-cost.

Clearly to bound~\eqref{eq:mse-mc} and the total computational cost,
we have to assume that we can control the FE error with respect to the
discretization parameter $h_\ell$ and that the computational cost to evaluate
the FE solution for a single sample $\rC \roundlr{{\rQ}_\ell(\by^{(m)})}$ is finite.
Combined with the third assumption given below,
we can express a bound for the total computational cost of the MLMC method
with respect to the target RMSE $\epsilon$ in the next section.

\begin{assumption}
    \label{assumption:mlmc}
    Suppose the approximation scheme to compute $\rQ_\ell$ satisfies
    \begin{align}
        \label{eq:assumption-alpha}
        \abs{\EE[{\rQ}_\ell - \rQ]} &\leq c_\alpha h_\ell^\alpha \\
        \label{eq:assumption-beta}
        \VV[{\rQ}_{\ell} - {\rQ}_{\ell - 1}] &\leq c_\beta h_\ell^\beta\\
        \label{eq:assumption-gamma}
        \rC \roundlr{{\rQ}_\ell(\by^{(m)})} &\leq c_\gamma h_\ell^{-\gamma}
    \end{align}
    with $\alpha, \beta, \gamma > 0$ and
    $c_\alpha, c_\beta, c_\gamma > 0$ independently on $\ell$.
\end{assumption}

\subsection{Introduction to Multi-level Monte Carlo methods}
\label{subsec:introduction-to-multi-level-monte-carlo}

The underlying idea of the MLMC method is to construct a model hierarchy for the imposed problem.
In FE applications, this can be done with nested meshes $\{\cM_\ell\}_{\ell=0}^L$ with decreasing mesh widths,
e.g.~$h_\ell =  2^{-\ell} h_0$, for the discretization on level $\ell=0,\dots,L$.
The goal is to reduce the number of evaluations of the model on
the finest level as much as possible and to minimize the overall
estimator variance. For a fixed finest level $L$, the expected value
of $\rQ_L$ can be written as a telescoping sum over the levels
\begin{equation}
    \label{eq:MLMC-telescope}
    \EE[{\rQ}_L] = \EE[{\rQ}_{0}] + \sum_{\ell=1}^L \EE[{\rQ}_{\ell} - {\rQ}_{\ell-1}]
    = \sum_{\ell=0}^L \EE[\rY_{\ell}]
    \,,
    \quad
    \rY_{0} \coloneqq {\rQ}_{0}
    \,,
    \quad
    \rY_{\ell} \coloneqq {\rQ}_{\ell} - {\rQ}_{\ell-1}
    \,.
\end{equation}
Each expected value of $\rY_\ell$ in the telescoping sum is now
estimated individually with a MC method, resulting in the MLMC
estimator
\begin{equation}
    \label{eq:mlmc-estimator-qoi}
    \widehat{{\rQ}}^{\text{MLMC}}_{\{M_\ell\}_{\ell=0}^L}
    = \sum_{\ell=0}^L \widehat{\rY}_{\ell}
    \,,
    \qquad
    \widehat{\rY}_{\ell}
    = \frac{1}{M_\ell} \sum_{m=1}^{M_\ell} \rY_\ell( \by^{(m)})\,,
\end{equation}
where $\{M_\ell\}_{\ell=0}^L$ denotes a sequence for the number of
samples on each level. It is important that every $\rY_\ell(
\by^{(m)}) = {\rQ}_{\ell}( \by^{(m)}) - {\rQ}_{\ell - 1}( \by^{(m)})$
uses the same sample $ \by^{(m)} \in \Xi$ for two different
meshes. Since all the expected values $\EE[\rY_\ell]$ are estimated
independently, the variance of the MLMC method can be quantified on
each level individually and with this, we obtain for the mean squared error
\begin{equation}
    \label{eq:mse-fem-mlmc}
    \mathrm{err}_{\mathrm{MSE}} \roundlr{\widehat{{\rQ}}^{\mathrm{MLMC}}_{\{M_\ell\}_{\ell=0}^L}}
    = \underbrace{\sum_{\ell=0}^L \frac{1}{M_\ell} \mathbb{V}[\rY_\ell]}_{\text{Estimator variance}}
    + \,\,\, {\underbrace{\EE[{\rQ}_L - {\rQ}]}_{\text{Bias}}}^2
    \,,
\end{equation}
cf. \cite{cliffe2011multilevel}.
Assuming we want to achieve a MSE of $\epsilon^2$, i.e., an RMSE
tolerance of $\epsilon$, we can reach this accuracy with $\theta \in
(0, 1)$, if
\begin{align*}
    \EE[{\rQ}_L - {\rQ}]^2 < (1 - \theta) \epsilon^2
    \quad \text{and} \quad
    \VV \squarelr{\widehat{{\rQ}}^{\text{MLMC}}_{\{M_\ell\}_{\ell=0}^L}} < \theta \epsilon^2.
\end{align*}
The parameter $\theta$ can thereby be used to tune the variance bias trade-off
in order to favor the minimization of one term over the other.
Numerical experiments have shown that for our particular
application $\theta = 0.5$ is a sufficient choice~\cite{baumgarten2023fully}.

Quantitatively, the computational cost of the method is given by
\begin{equation}
    \label{eq:mlmc-cost}
    \rC \roundlr{\widehat{{\rQ}}^{\text{MLMC}}_{\{M_\ell\}_{\ell=0}^L}}
    = \sum_{\ell=0}^L \sum_{m=1}^{M_\ell} \rC_{\ell}(\by^{(m)})
    = \sum_{\ell=0}^L M_\ell \widehat{\rC}_{\ell}
    \,,
    \quad
    \widehat{\rC}_{\ell}
    =
    \frac{1}{M_\ell}
    \sum_{m=1}^{M_\ell} \rC_{\ell}(\by^{(m)})
    \,,
\end{equation}
where $\widehat{\rC}_{\ell}$ is the sample mean of the cost. We now want to find the
optimal sequence of samples $\{M_\ell\}_{\ell=0}^L$, such that the
estimator cost is minimized while achieving an MSE tolerance of
$\epsilon^2$, cf.~\cite{giles2015multilevel}. By presetting the MSE tolerance,
we can also deduce from~\eqref{eq:assumption-alpha} the highest level~$L$ since the
estimator bias has to be smaller than $(1 - \theta)
\epsilon^2$. Thereby, it is sufficient to minimize the estimator cost
while achieving an estimator variance of $\theta \epsilon^2$, i.e., we
search for $\{M_\ell\}_{\ell=0}^L$ solving
\begin{equation}
    \label{eq:cost-minimization}
    \min_{\{M_\ell\}_{\ell=0}^L} \quad \rC \roundlr{\widehat{{\rQ}}^{\text{MLMC}}_{\{M_\ell\}_{\ell=0}^L}}
    \quad \text{s.t.} \quad
    \sum_{\ell=0}^L M_\ell^{-1} \VV[\rY_\ell] = \theta {\epsilon^2}.
\end{equation}
By treating each $M_\ell$ as a continuous variable, the solution to this optimization problem is given by
\begin{equation}
    \label{eq:optimal-Ml}
    M_\ell^{\mathrm{opt}} = \ceil { \roundlr{{\sqrt{\theta} \epsilon}}^{-2}
        \sqrt{\frac{\VV[\rY_\ell]}{\widehat{\rC}_\ell}}
        \roundlr{\sum_{\ell'=0}^L \sqrt{\VV[\rY_{\ell'}] \widehat{\rC}_{\ell'}}}} \quad
    \text{for} \quad \ell = 0, \dots, L.
\end{equation}
We lastly restate the $\epsilon$-cost theorem of the MLMC method
in the form given in~\cite{cliffe2011multilevel, giles2008improved}:

\begin{theorem}
    [Bounded $\epsilon$-cost of the MLMC method]
    \label{theorem:mlmc-epsilon-theorem}
    Suppose assumption~\ref{assumption:mlmc} is fulfilled with some
    positive rates $\alpha, \beta, \gamma > 0$ with $\alpha \geq
    \tfrac{1}{2} \min \set{\beta, \gamma}$ and $c_\alpha, c_\beta, c_\gamma > 0$
    independent of~$h_\ell$. Then for any $0 < \epsilon < \re^{-1}$, there
    exists a maximum level $L \in \NN$ and a sequence of samples
    $\{M_\ell\}_{\ell=0}^L$ such that
    \begin{align*}
        {\mathrm{err}_{\mathrm{RMSE}}} \roundlr{\widehat{{\rQ}}^{\mathrm{MLMC}}_{\{M_\ell\}_{\ell=0}^L}} < \epsilon
        \qquad \text{with} \qquad
        \rC_\epsilon \roundlr{\widehat{{\rQ}}^{\mathrm{MLMC}}_{\{M_\ell\}_{\ell=0}^L}} {\lesssim}
        \begin{cases}
            \epsilon^{-2} & \beta > \gamma\,, \\
            \epsilon^{-2}\log(\epsilon)^2 & \beta = \gamma\,, \\
            \epsilon^{-2-(\gamma - \beta)/\alpha} & \beta < \gamma\,.
        \end{cases}
    \end{align*}
    with $\lesssim$ denoting that the left hand side obeys an upper bound given
    by the right hand side up to some hidden constant.
\end{theorem}

\subsection{Implementation Techniques}
\label{subsec:implementation-techniques}

Implementing the MLMC method with on-the-fly estimation
of~\cref{assumption:mlmc} is well described
in~\cite{giles2015multilevel} and shortly recalled to commit to our
notation. This includes a way to find the optimal sequence
$\{M_\ell\}_{\ell=0}^L$ and the highest level $L$ during runtime.
Choosing $M_\ell^{\mathrm{opt}}$ by~\eqref{eq:optimal-Ml} 
requires estimates for the variance $\VV[\rY_\ell]$ and the sample mean of the cost
$\widehat{\rC}_\ell$, which are unknown a priori. The idea is to perform
the MLMC method with an initial sequence
$\{M_{\ell}^{\text{init}}\}_{\ell=0}^{L_{\text{init}}}$ to get first
estimates for the sample mean of the cost $\widehat{\rC}_{\ell}$ and
the variance $\VV[\rY_\ell]$ by the sample variance estimator
\begin{align}
  \label{eq:sample-variance}
  s_{\rY_\ell}^2 = \frac{1}{M_{\ell} - 1} S_{\rY_\ell, 2} \quad \text{with}
  \quad S_{\rY_\ell, 2} = \sum_{m=1}^{M_\ell} \big(\rY_\ell^{(m)} - \widehat{\rY}_\ell\big)^2
  \,.
\end{align}
From these initial estimates on, the MLMC method is executed until the target
RMSE $\epsilon$ is reached by continuously updating the sample statistics.
Hence, the amount of samples accumulates and the required
amount is given by ${\vartriangle} M_{\ell} \coloneqq \max\{\widehat{M}_{\ell}^{\text{opt}} - M_{\ell}, \, 0\}$ with
the optimal sample amount based on the estimates
\begin{equation}
    \label{eq:optimal-Ml-estimated}
    \widehat{M}_{\ell}^{\text{opt}} = \ceil {
        \roundlr{\sqrt{\theta} \epsilon}^{-2}
        \sqrt{\frac{s^2_{\rY_{\ell}}}{\widehat{\rC}_{\ell}}}
        \roundlr{\sum_{\ell'=0}^L \sqrt{s^2_{\rY_{\ell'}} \widehat{\rC}_{\ell'}}}} \quad
    \text{for} \quad \ell = 0, \dots, L\,.
\end{equation}
With assumption~\eqref{eq:assumption-alpha} and the geometric sum, an estimate for the bias is given by
\begin{equation}
    \label{eq:disc-estimate}
    \widehat{\err}_{\text{disc}} = \max \set{\frac{\widehat{\rY}_\ell}{2^{\widehat{\alpha}} - 1}
    2^{-{\widehat{\alpha}} (L - \ell)} \colon \ell = 1, \dots, L},
\end{equation}
which also incorporates lower levels for robustness of the estimate, cf.~\cite{giles2015multilevel},
and uses an approximation for ${{\alpha}}$ by fitting the data
$\{\widehat{\rY}_\ell\}_{\ell=0}^L$ to assumption~\eqref{eq:assumption-alpha}
\begin{equation}
  \min_{(\widehat{\alpha}, \, \widehat{c}_\alpha)} \quad \sum_{\ell=1}^L
  \Big(\log_2 \widehat{\rY}_\ell + \widehat{\alpha} \ell - \widehat{c}_\alpha \Big)^2
  \,.
    \label{eq:alpha-fit}
\end{equation}
This gives the estimates $\widehat{\alpha}$ and $\widehat{c}_\alpha$
for the rate $\alpha$ and the constant $c_\alpha$ which
can be done in a similar way for $(\gamma, c_\gamma)$~and~$(\beta, c_\beta)$.  The estimator
variance of the MLMC method is approximated with the sample variance by 
\begin{equation*}
  \widehat{\err}_{\text{input}} = 
  \sum_{\ell=0}^L \frac{1}{M_{\ell}} s^2_{\rY_\ell}
\end{equation*}
and, finally, an estimate for the MSE can be computed with
\begin{equation*}
  \widehat{\err}_{\text{MSE}} = \widehat{\err}_{\text{input}} +
  \widehat{\err}_{\text{disc}}^2
  \,.
\end{equation*}
The above techniques are further refined by the CMLMC
method~\cite{collier2015continuation}, also applied
in~\cite{litvinenko2019computation, pisaroni2017continuation}. The
central idea is to create a sequence $\epsilon_\ttti$
with $\epsilon_{\ttti} \rightarrow \epsilon_{\min} > 0$ and by doing
so, the above estimates and the required number of samples
$\{M_{\ttti,\ell}\}_{\ell=0}^{L_\ttti}$ are continuously updated for each new
tolerance. This idea is now extended with an additional constraint on
the total cost, i.e.,~a budget, giving the budgeted MLMC method.

\subsection{Budgeted Multi-level Monte Carlo}\label{subsec:budgeted-multi-level-monte-carlo}

In large HPC systems, the workload managers require to reserve
computational budget $\rB \coloneqq \abs{\cP} \cdot \rT_{\mathrm{\rB}}$
to initiate a job on the cluster. Therefore, we define an execution
of a method on a computing cluster as \textit{feasible},
if the method is capable of fully utilizing the computational budget without exceeding it.
In light of this, we adapt the perspective of~\cref{theorem:mlmc-epsilon-theorem}.
The final implementation will contain certain parts of the algorithm which can be parallelized,
while others cannot. We denote $\lambda_\rp \in [0, 1]$ as the parallelization constant,
representing the portion of the implemented algorithm that is executed in parallel on $\abs{\cP}$ units.

\begin{proposition}[Convergence of a parallelized BMLMC method]
    \label{proposition:bmlmc-proposition}
    For a feasible execution of the budgeted method by a parallel
    implementation, the estimate for the error splits up into two parts
    \begin{equation}
        \label{eq:bmlmc-theorem}
        \epsilon
        \lesssim \underbrace{(1 - \lambda_\rp) \cdot \rT_{\mathrm{\rB}}^{-\delta}}_{\eqqcolon \epsilon_\rs}
        + \underbrace{\lambda_\rp (\abs{\cP} \cdot \rT_{\mathrm{\rB}})^{-\delta}}_{\eqqcolon \epsilon_\rp}
        \quad \text{with} \quad \delta =
        \begin{cases}
            \hspace{0.6cm}
            \frac{1}{2} & \beta > \gamma \\
            \frac{\alpha}{2 \alpha + (\gamma - \beta)} & \beta < \gamma
        \end{cases}
    \end{equation}
    depending on the parallization constant $\lambda_\rp \in [0, 1]$.
\end{proposition}

For the proof we refer to~\Cref{subsec:parallelization-techniques}, here we only
comment on the case $\lambda_\rp =1$ when a perfectly parallel
implementation can be realized. Then, in case of a feasible run,
the final $\epsilon$-cost equals the budget, and the estimate
\begin{align}
    \label{eq:lambda1}
    \epsilon \lesssim \rB^{-\delta} =
    \begin{cases}
        \rB^{-1/2} & \beta > \gamma \\
        \rB^{-\alpha/(2 \alpha + (\gamma - \beta))} & \beta < \gamma
    \end{cases}
\end{align}
simply follows by inverting~\cref{theorem:mlmc-epsilon-theorem}.
Similar results can be found e.g.~in~\cite{kumar2018multigrid,
    mishra2016multi} where the word \textit{work} was used instead of
budget. The case $\beta = \gamma$ is neglected for the sake of a leaner representation
and since it has no practical relevance if the cost is measured in units of CPU seconds.

This result does not tell us how to utilize the budget. Thereby, the
new algorithmic challenge is to find the best way to invest $\rB$,
such that we minimize the error.
Formally, this is expressed by a knapsack problem:

\begin{problem}[MLMC Knapsack]
    \label{problem:mlmc-knapsack}
    Find $L$ and $\{M_\ell\}_{\ell=0}^L$, such that the MSE is minimized while staying within
    the cost budget $\rB$,~i.e.,
    \begin{subequations}
        \begin{align}
            \label{eq:knapsack-mse}
            \min_{(L, \{M_\ell\}_{\ell=0}^L)} \quad &\err_{\mathrm{MSE}}
            = \sum_{\ell=0}^L \frac{1}{M_\ell} \VV[\rY_\ell] + \big(\EE[{\rQ}_L - {\rQ}]\big)^2 \\
            \text{s.t.} \quad
            &\sum_{\ell = 0}^L \sum_{m=1}^{M_{\ell}} \rC_{\ell}(\by^{(m)}) \leq \rB
            \label{eq:knapsack-restriction}.
        \end{align}
    \end{subequations}
\end{problem}

We remark the change of perspective
to~\eqref{eq:cost-minimization}. The objective now is to minimize the
complete MSE including the bias, while the constraint is given by the
computational budget. Knapsack problems, in general, are combinatorial
optimization problems and often arise while searching for the optimal
allocation of resources, e.g., in manufacturing, in computer networks or
in financial models. These integer optimization problems are NP-hard
and require effective algorithms often designed with dynamic
programming (DP) techniques~\cite{bretthauer2002nonlinear,
    caccetta2001computational}. The key idea of DP is to split up the
initial problem into \textit{overlapping subproblems} and solve them
recursively with some \textit{optimal policy}, while reusing
    \textit{memoized} results stored in a suitable data structure.

Now, the goal is to derive an algorithm satisfying \textit{Bellman's
optimality condition} of DP, i.e., an algorithm finding a sequence
of optimal actions, such that at each state an objective value is
maximized. To do so, we have to construct a \textit{reward/pay-off} function
which describes in each state the reward/pay-off for the maximization,
if a certain action is taken.
For the construction of a suitable reward function, we recognize
that we do not know the exact quantities in~\eqref{eq:knapsack-mse}
and~\eqref{eq:knapsack-restriction}.
The following optimization problem is
still NP-hard but at least exclusively contains computable quantities.

\begin{problem}[Approximated MLMC Knapsack]
    \label{problem:approximated-mlmc-knapsack}
    Find $L$ and the sequence $\{M_\ell\}_{\ell=0}^L$, such that the estimated MSE is minimized,
    while staying within the cost budget $\rB$, i.e.,
    \begin{subequations}
        \begin{align}
            \label{eq:knapsack-estimated-mse}
            \min_{(L, \{M_\ell\}_{\ell=0}^L)} \quad &\widehat{\err}_{\mathrm{MSE}}
            = \sum_{\ell=0}^L \frac{1}{M_{\ell}} s^2_{\rY_\ell}
            + \widehat{\err}_{\mathrm{disc}}^2 \\
            \text{s.t.} \quad
            &\sum_{\ell = 0}^L M_{\ell} \widehat{\rC}_{\ell} \leq \rB
            \,,
            \label{eq:knapsack-estimated-restriction}
        \end{align}
    \end{subequations}
    where $\widehat{\err}_{\mathrm{disc}}$ is computed with~\eqref{eq:disc-estimate}
    and the sample variance $s_{\rY_\ell}^2$ by \eqref{eq:sample-variance}.
\end{problem}
The idea is to identify the MSE as the value we try to optimize and to
split up the initial problem into several estimation rounds with a
decreasing sequence of tolerances~$\epsilon_\ttti$
as in the CMLMC method~\cite{collier2015continuation}. This
creates subsequent optimization problems where each solution yields
some reward to the total optimization. To this end, we equip all
quantities of~\Cref{subsec:implementation-techniques} with an
index $\ttti$. In the following, we motivate the existence of a
\textit{pay-off/reward} function depending on the chosen action
$\epsilon_{\ttti}$, and the current state of the simulation, i.e., the
collected data up to $\ttti - 1$.

Suppose we are in estimation round $\ttti$ and $M_{0,
\ell}^{\text{init}}$ is given. For $\ell=0,\dots,L_{\ttti}$, we separate
\begin{align*}
    M_{\ttti, \ell} &= \underbrace{M_{\ttti - 1, \ell}}_{\text{available data}}
    + \underbrace{{\vartriangle} M_{\ttti, \ell}}_{\text{optimal choice}}
    \quad \text{with} \quad {\vartriangle} M_{0, \ell}
    \coloneqq M_{0, \ell}^{\text{init}} \quad \text{and} \quad M_{-1, \ell} \coloneqq 0
\end{align*}
such that ${\vartriangle} M_{\ttti, \ell}$
is computed using~\eqref{eq:optimal-Ml-estimated}, hence after the estimation
round we have $M_{\ttti, \ell} = \widehat{M}_{\ttti,
\ell}^{\text{opt}}$. Thereby, we can express the amount of samples
based on the currently available data and some \textit{optimal
policy}, i.e., ${\vartriangle} M_{\ttti, \ell}$ is chosen such that
the cost is minimized and a target MSE tolerance of $\epsilon_\ttti^2$
is reached. With ${\vartriangle} M_{\ttti, \ell}$ and accumulative
update formulas for the sample mean and sample
variance~\cite{pebay2016numerically}, we further separate
\begin{align*}
    \widehat{\rC}_{\ttti, \ell} &= \widehat{\rC}_{\ttti - 1, \ell}
    + \frac{{\vartriangle} M_{\ttti, \ell}}{\widehat{M}_{\ttti, \ell}^{\text{opt}}}
    ({\vartriangle} \widehat{\rC}_{\ttti, \ell} - \widehat{\rC}_{\ttti - 1, \ell})
    \quad \text{with} \quad
    {{\vartriangle} \widehat{\rC}_{\ttti, \ell} \coloneqq \frac{1}{{\vartriangle} M_{\ttti, \ell}}
    \sum\limits_{m = M_{\ttti - 1, \ell} + 1}^{\widehat{M}_{\ttti, \ell}^{\text{opt}}}  \rC_\ell^{(m)}} \\
    \widehat{\rY}_{\ttti, \ell} &= \widehat{\rY}_{\ttti - 1, \ell}
    + \frac{{\vartriangle} M_{\ttti, \ell}}{\widehat{M}_{\ttti, \ell}^{\text{opt}}}
    ({\vartriangle} \widehat{\rY}_{\ttti, \ell} - \widehat{\rY}_{\ttti - 1, \ell})
    \quad \text{with} \quad
    {{\vartriangle} \widehat{\rY}_{\ttti, \ell} \coloneqq \frac{1}{{\vartriangle} M_{\ttti, \ell}}
    \sum\limits_{m = M_{\ttti - 1, \ell} + 1}^{\widehat{M}_{\ttti, \ell}^{\text{opt}}}  \rY_\ell^{(m)}}
\end{align*}
and likewise for
\begin{align*}
    s_{\rY_{\ttti, \ell}}^2 =
    \frac{S_{\rY_{2, \ttti, \ell}}}{\widehat{M}_{\ttti, \ell}^{\text{opt}} - 1} ,
\end{align*}
we separate
\begin{align*}
    S_{\rY_{2, \ttti, \ell}} =  S_{\rY_{2, \ttti - 1, \ell}} + {\vartriangle} S_{\rY_{2, \ttti, \ell}}
    + \frac{M_{\ttti - 1, \ell} {\vartriangle} M_{\ttti, \ell}}{\widehat{M}_{\ttti, \ell}^{\text{opt}}}
    \big(\widehat{\rY}_{\ttti - 1, \ell} - {\vartriangle} \widehat{\rY}_{\ttti, \ell}\big)^2\,,
\end{align*}
where
\begin{align*}
{\vartriangle}
    S_{\rY_{2, \ttti, \ell}} = \sum\limits_{m = M_{\ttti - 1 , \ell} + 1}^{\widehat{M}^{\text{opt}}_{\ttti, \ell}} \roundlr{\rY_\ell^{(m)} - {\vartriangle} \widehat{\rY}_{\ttti, \ell}}^2.
\end{align*}
As $\widehat{\alpha}_{\ttti}$ is computed with~\eqref{eq:alpha-fit}, i.e.,
a fit to the available data, we can express~\eqref{eq:knapsack-estimated-mse}
and~\eqref{eq:knapsack-estimated-restriction} for a particular estimation round $\ttti$
as a nonlinear function of preexisting data (the \textit{state}) and
the \textit{optimal policy}~\eqref{eq:optimal-Ml-estimated}.
This motivates a function for the \textit{pay-off} purely determined by the state and the optimal policy.
We denote this function by ${\vartriangle} \err_{\text{MSE}}(\texttt{data}_{\ttti - 1}, {\vartriangle} \texttt{data}_{\ttti})$
which represents the error reduction in one estimation round,
if $\{{\vartriangle} M_{\ttti, \ell}\}_{\ell=0}^{L_\ttti}$ additional samples are computed. We use the notation
\begin{align*}
    \texttt{data}_{\ttti}
    = \set{ \texttt{data}_{\ttti, \ell}}_{\ell=0}^{L_{\ttti}}
    &= \set{\set{ M_{\ttti, \ell},
        \widehat{\rQ}_{\ttti, \ell},
        \widehat{\rC}_{\ttti, \ell},  \widehat{\rY}_{\ttti, \ell},
        S_{\rQ_{2, \ttti, \ell}}, S_{\rY_{2, \ttti, \ell}}, \dots}}_{\ell=0}^{L_{\ttti}}
    \,, \\
    {\vartriangle} \texttt{data}_{\ttti}
    = \set{{\vartriangle} \texttt{data}_{\ttti, \ell}}_{\ell=0}^{L_{\ttti}}
    &= \set{\set{{\vartriangle} M_{\ttti, \ell},
            {\vartriangle} \widehat{\rQ}_{\ttti, \ell},
            {\vartriangle} \widehat{\rC}_{\ttti, \ell},
            {\vartriangle} \widehat{\rY}_{\ttti, \ell},
            {\vartriangle} S_{\rQ_{2, \ttti, \ell}},
            {\vartriangle} S_{\rY_{2, \ttti, \ell}}, \dots}}_{\ell=0}^{L_{\ttti}}
\end{align*}
to collect all needed quantities in one object.
We further define $\rB_\ttti$ as the left-over budget in round $\ttti$
and denote with $\rB_0 \coloneqq \rB$ the initially imposed
budget. By~\eqref{eq:optimal-Ml-estimated}, we see that the amount of
samples is guided by $\epsilon_{\ttti}$.
With $\epsilon_{\ttti}$ as the chosen \textit{action} and
$\eta \in (0, 1)$ as reduction factor determining how fast $\epsilon_\ttti$ decays
and with the cost prediction
$\widehat{\rC}_{\ttti} = \sum_{\ell=0}^{L_\ttti} {\vartriangle} M_{\ttti, \ell} \widehat{\rC}_{\ttti - 1, \ell}$,
the \textit{Bellman equation} for finding the solution to
\cref{problem:approximated-mlmc-knapsack} can be expressed by
\begin{align*}
    \widehat{\err}_{\text{MSE}}^{\text{final}}(\rB_0, \{M_{0,\ell}^{\text{init}}\}_{\ell=0}^{L_0}) =
    \widehat{\err}_{\text{MSE}}^{\text{init}}
    - \widehat{\err}_{\text{MSE}} \roundlr{\rB_{0} - \textstyle
    \sum\limits_{\ell=0}^{L_0} \rC_{0,\ell}, \, \eta \cdot \widehat{\err}_{\text{MSE}}^{\text{init}}}
\end{align*}
with the recursive function
\begin{align}
    \label{eq:bmlmc-bellman}
    \widehat{\err}_{\text{MSE}}(\rB_\ttti, \epsilon_\ttti) =
    \max_{\substack{\set{{\vartriangle} M_{\ttti, \ell}}_{\ell=0}^{L_\ttti} \\
    \text{s.t. } \widehat{\rC}_{\ttti} < \rB_\ttti}}
    \Big\{&{\vartriangle} \err_{\text{MSE}}(\texttt{data}_{\ttti - 1}, {\vartriangle} \texttt{data}_{\ttti}) \\[-6mm]
    &\quad + \, \widehat{\err}_{\text{MSE}}
    \Big(\rB_{\ttti} - \textstyle \sum\limits_{\ell=0}^{L_\ttti} \rC_{\ttti, \ell},
    \, \eta \cdot \epsilon_{\ttti}\Big) \Big\}.
    \nonumber
\end{align}
So far, we have not discussed the minimization of the bias yet.
If in~\eqref{eq:knapsack-estimated-mse} the bias becomes larger than
$\sqrt{1 - \theta} \, \epsilon_\ttti$
and if we have enough budget left, i.e., $\widehat{\rC}_\ttti < \rB_\ttti$,
we draw additional samples on level $L_{\ttti}+1$ and stop the optimization otherwise.

\begin{figure}[H]
    \centering
    \includegraphics[width=0.4\textwidth]{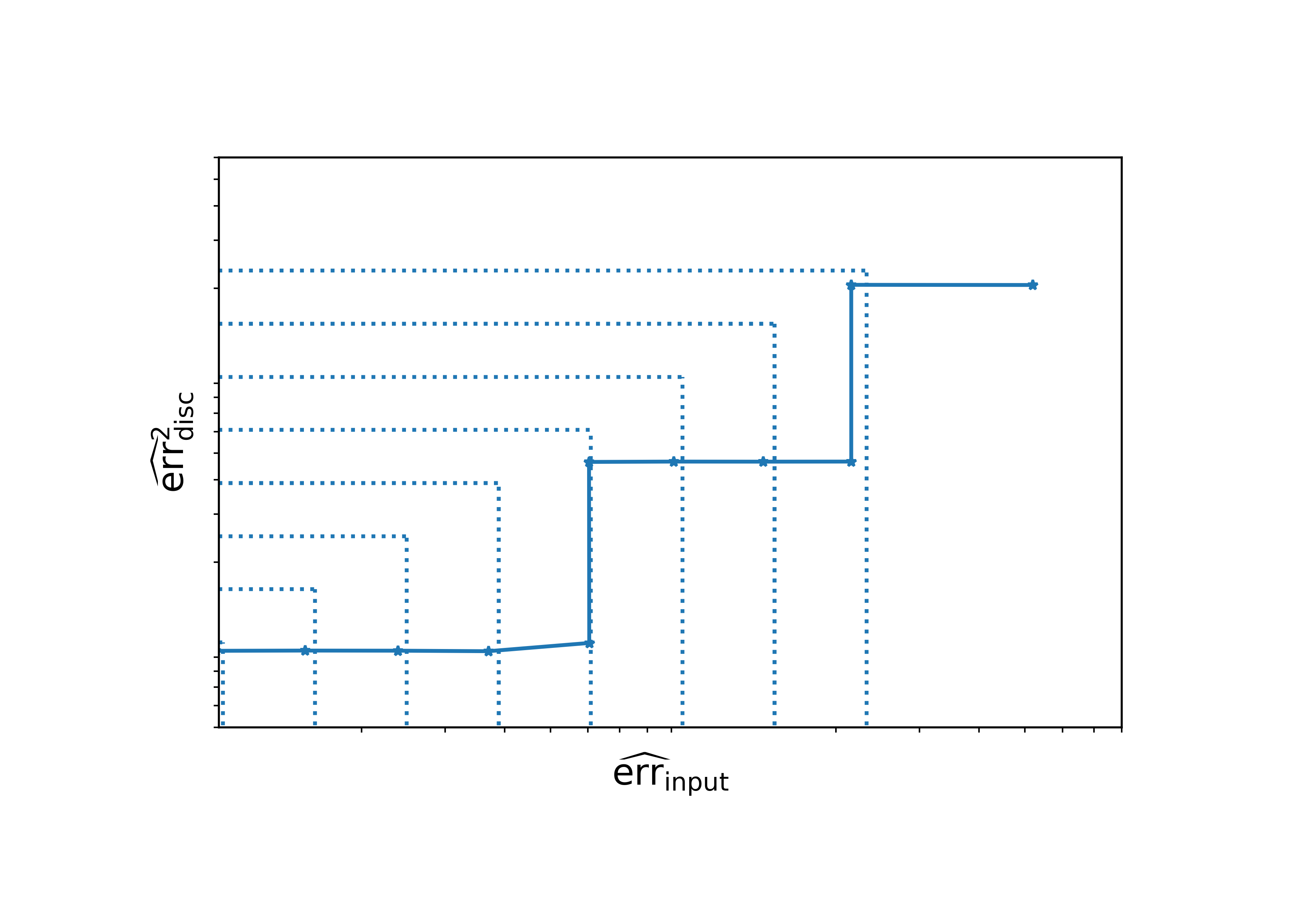}
    \label{fig:mlmc-elliptic-l2-on-horeka-00002-11-1}
    \caption{Illustration of an example execution
    of~\cref{alg:bmlmc}. The data point in the upper right corner
    represents the first estimate of
        $\widehat{\err}_{\mathrm{disc}}^2$ and
        $\widehat{\err}_{\mathrm{input}}$ based on the initial sequence
        $\{M_{0, \ell}^{\mathrm{init}}\}_{\ell=0}^{L_0}$.
        The doted
        squares represent the target MSEs $\epsilon_{\ttti}^2$ in each
        estimation round and illustrate the overlapping
        subproblems. Hence, \cref{alg:bmlmc} solves in each estimation
        round the optimization problem~\eqref{eq:cost-minimization} or
        adapts the highest level $L_{\ttti}$. The joint effort of each
        estimation round finally solves
        \cref{problem:approximated-mlmc-knapsack}.}
\end{figure}

To conclude, function~\eqref{eq:bmlmc-bellman} is the expression of
\textit{Bellman's optimality condition} applied to
\cref{problem:approximated-mlmc-knapsack}, i.e., the subsequent
minimization of the MSE under consideration of the cost budget.
This subsequent minimization is also illustrated
in~\cref{fig:mlmc-elliptic-l2-on-horeka-00002-11-1},
where each dotted square represents one estimation round.

\begin{figure}[H]
    \centering
    \includegraphics[width=0.4\textwidth]{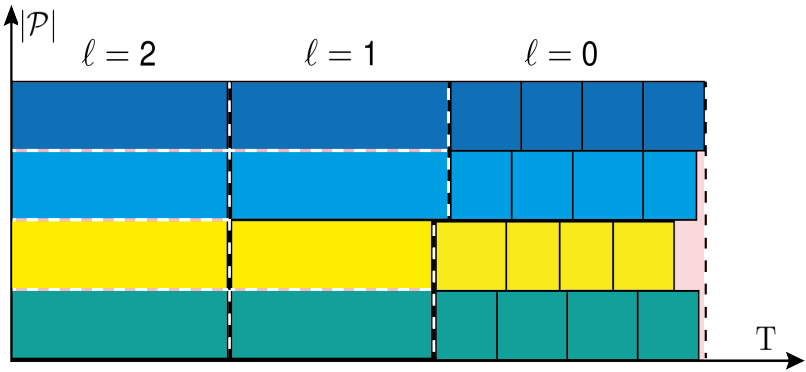}
    \caption{Load balancing of a single exemplary estimation round on three levels with $M_{0} = 16$, $M_{1} = 2$ and $M_{2} = 1$. The sample on level two is processed on all four processes at first. Next, two samples on level one follow where each one is computed on two processes. Lastly sixteen samples on level zero are computed, where each processor handles four individual samples. The light red areas correspond to parallelization losses due to either idling processes (red areas) or due to communication losses (dashed red lines, corresponding to communication across subdomains). The remaining colors repent one of four processing units. The black dashed lines represent transitions between levels, the red dashed lines represent parallelization losses in the FEM system.}
    \label{fig:load-balancing}
\end{figure}

\begin{remark}
    The actual implementation in C++~\cite{wieners2023mpp315} is not done
    with a recursive function but in an equivalent formulation with a
    while-loop. This is also often called a \textit{bottom up}
    implementation which has the advantage over the recursive
    implementation (\textit{top down}) to avoid an increased memory
    consumption on the stack. However, the recursive formulation is easier
    to derive mathematically. We further remark the inverted level loop
    which has benefits for the load distribution~\cite{badia2021massively}
    as illustrated in~\cref{fig:load-balancing}.
\end{remark}

Lastly, we present in~\Cref{alg:bmlmc} the final BMLMC method as a recursive implementation.
For a detailed explanation of the subroutines \texttt{Welford} and \texttt{MS-FEM},
we refer to the upcoming \Cref{subsec:parallelization-techniques} and to \Cref{sec:multi-sample-finite-element-method}.

\begin{algorithm}
    \caption{Budgeted Multi-Level Monte Carlo}
    \label{alg:bmlmc}
    Set the initial sample sequence $\{M_{0, \ell}^{\text{init}}\}_{\ell=0}^{{L_0}}$
    in $\texttt{data}_{-1}$, a cost budget $\abs{\cP} \cdot \rT_{\rB} = \rB > 0$, a splitting parameter $\theta \in (0, 1)$ and the reduction factor $\eta \in (0, 1)$.
    \begin{align*}
        &\texttt{data} = \set{\ttti \mapsto \set{\err_{\ttti}, \{M_{\ttti, \ell}\}_{\ell=0}^{L_{\ttti}},
            \{\widehat{\rQ}_{\ttti, \ell}\}_{\ell=0}^{L_{\ttti}},
            \{\widehat{\rC}_{\ttti, \ell}\}_{\ell=0}^{L_{\ttti}},
            \{\widehat{\rY}_{\ttti, \ell}\}_{\ell=0}^{L_{\ttti}}, \dots}} \\[2mm]
        &\texttt{function } \texttt{BMLMC}(\rB_{0}, \set{M_{0, \ell}^{\text{init}}}_{\ell=0}^{L_0}) \colon \\
        &\quad
        \begin{cases}
            \texttt{for } \ell = {L_0}, \dots, 0 \colon
            &\hspace{0.3cm} {\vartriangle} \texttt{data}_{0, \ell} \leftarrow \texttt{MS-FEM}(M_{0, \ell}^{\text{init}}, \cP) \\
            \texttt{data}_{0} \leftarrow \texttt{Welford}(\texttt{data}_{-1}, {\vartriangle} \texttt{data}_{0})
            &\hspace{0.3cm} \texttt{return } \texttt{BMLMC}(\rB_{0} - \sum_{\ell=0}^{L_{0}} \rC_\ell, \eta \cdot \err_{0})
        \end{cases} \\[2mm]
        &\texttt{function } \texttt{BMLMC}(\rB_{\ttti}, \epsilon_{\ttti}) \colon \\
        &\quad
        \begin{cases}
            \texttt{if } \rB_{\ttti} \approx 0 \colon
            &\hspace{0.2cm} \texttt{return } \err_{\ttti - 1} \\
            \texttt{if } \widehat{\err}_{\texttt{disc}}(\texttt{data}_{\ttti - 1}) \geq \sqrt{1 - \theta} \epsilon_{\ttti} \colon
            &\hspace{0.2cm} {L_{\ttti}} \leftarrow {L_{\ttti}} + 1 \\
            \texttt{if } \widehat{\err}_{\texttt{input}}(\texttt{data}_{\ttti - 1}) \geq \theta \epsilon_{\ttti}^2 \colon
            &\hspace{0.2cm} \widehat{M}_{\ttti, \ell}^{\text{opt}} \leftarrow~\eqref{eq:optimal-Ml-estimated} \\
            \texttt{for } \ell = {L_{\ttti}}, \dots, 0 \colon
            &\hspace{0.2cm}{\vartriangle} M_{\ttti, \ell} \leftarrow \max \set{\widehat{M}_{\ttti, \ell}^{\text{opt}} - M_{\ttti-1,\ell}, \, 0} \\
            \widehat{\rC}_{\ttti} \leftarrow \sum_{\ell=0}^{L_{\ttti}} {\vartriangle} M_{\ttti, \ell} \widehat{\rC}_{\ttti - 1, \ell} \\
            \texttt{if } \widehat{\rC}_{\ttti} = 0 \colon
            &\hspace{0.2cm} \texttt{return } \texttt{BMLMC}(\rB_{\ttti}, \eta \cdot \epsilon_{\ttti}) \\
            \texttt{if } \widehat{\rC}_{\ttti} > \rB_{\ttti} \colon
            &\hspace{0.2cm} \texttt{return } \texttt{BMLMC}(\rB_{\ttti}, 0.5 \cdot (\epsilon_{\ttti} + \epsilon_{\ttti - 1})) \\
            \texttt{for } \ell = {L_{\ttti}}, \dots, 0 \colon
            &\hspace{0.2cm} {\vartriangle} \texttt{data}_{\ttti, \ell} \leftarrow \texttt{MS-FEM}({\vartriangle} M_{\ttti, \ell}, \cP) \\
            \texttt{data}_{\ttti} \leftarrow \texttt{Welford}(\texttt{data}_{\ttti - 1}, {\vartriangle} \texttt{data}_{\ttti})
            &\hspace{0.2cm} \texttt{return } \texttt{BMLMC}(\rB_{\ttti} - \sum_{\ell=0}^{L_{\ttti}} \rC_\ell, \epsilon_{\ttti})
        \end{cases}
    \end{align*}
\end{algorithm}

\begin{remark}
    By considering \cref{problem:approximated-mlmc-knapsack}, we chose to
    \textit{discretize first and optimize then}. The downside of this
    approach is that if $s_{\rY_\ell}^2$, $\widehat{\rC}_\ell$,
    $\widehat{\rY}_\ell$ and $\widehat{\alpha}$ are inaccurate, the
    optimization delivers poor results, too. However, by using dynamic
    programming we actually \textit{discretize, optimize, discretize,
        optimize, \dots} until the budget is exhausted. Hence, the risk of
    optimizing for the wrong objective based on inaccurate data is reduced
    as the simulation runs.
\end{remark}

\subsection{Parallelization Techniques}\label{subsec:parallelization-techniques}

\cref{alg:bmlmc} is designed as a state machine.
However, due to the significant computational load involved in solving
\cref{problem:approximated-mlmc-knapsack}, it becomes necessary to distribute the workload
across multiple nodes or processing units. As a result, \cref{alg:bmlmc} needs to be
transformed into a \textit{distributed state machine}, i.e., a computational method
on interconnected nodes or processing units that synchronizes and maintains a shared
state. In this context, the shared state refers to the
data containing the estimated errors and sample statistics in the very first line of the algorithm, while the set of
processing units $\cP$ is responsible for dividing the work and minimizing $\widehat{\err}_{\text{MSE}}$.

\begin{algorithm}
    \caption{Welford's Online Update Method}
    \label{alg:welfords-weighted-online-update}
    Compute for
    $\big\{\{\texttt{data}_{\ttti, \ell}\}_{\ell=0}^{L_{\ttti}}\big\}_{\ttti = A,B} =
    \big\{\{M_{\ttti, \ell}, \widehat{\rQ}_{\ttti, \ell}, S_{\rQ_{2, \ttti, \ell}}, \dots\}_{\ell=0}^{L_{\ttti}}\big\}_{\ttti = A,B}$
    \vspace{-2mm}
    \begin{align*}
        &\texttt{function Welford} (\texttt{data}_A, \texttt{data}_B)) \colon \\
        &\begin{cases}
             M_{AB, \ell} &\quad \leftarrow \quad M_{B, \ell} + M_{A, \ell} \\
             \delta_{AB, \ell} &\quad \leftarrow \quad \widehat{\rQ}_{B, \ell} - \widehat{\rQ}_{A, \ell} \\
             \widehat{\rQ}_{AB, \ell} &\quad \leftarrow
             \quad \widehat{\rQ}_{A, \ell} + \tfrac{M_{B, \ell}}{M_{AB, \ell}} \delta_{AB, \ell} \\
             S_{\rQ_{2, AB, \ell}} &\quad \leftarrow
             \quad S_{\rQ_{2, A, \ell}} + S_{\rQ_{2, B, \ell}} + \tfrac{M_{A, \ell} M_{B, \ell}}{M_{AB, \ell}} \delta_{AB, \ell}^2 \\
             s_{\rQ_{AB, \ell}}^2 &\quad \leftarrow
             \quad (M_{AB, \ell} - 1)^{-1} S_{\rQ_{2, AB, \ell}} \\[1mm]
             \texttt{return } &\{M_{AB, \ell}, \widehat{\rQ}_{AB, \ell}, \, S_{\rQ_{2, AB, \ell}}, \dots\}_{\ell=0}^{L_{\ttti}}
        \end{cases}
    \end{align*}
\end{algorithm}

To achieve this, we discuss the functionality of \Cref{alg:welfords-weighted-online-update} representing
the subroutine \texttt{Welford} in \Cref{alg:bmlmc}. We use concepts introduced
by~\cite{chan1983algorithms, welford1962note},
which were further expanded in~\cite{pebay2016numerically},
to stably compute sample statistics in an incremental and parallel manner.

\begin{figure}[H]
    \centering
    \includegraphics[width=0.75\textwidth]{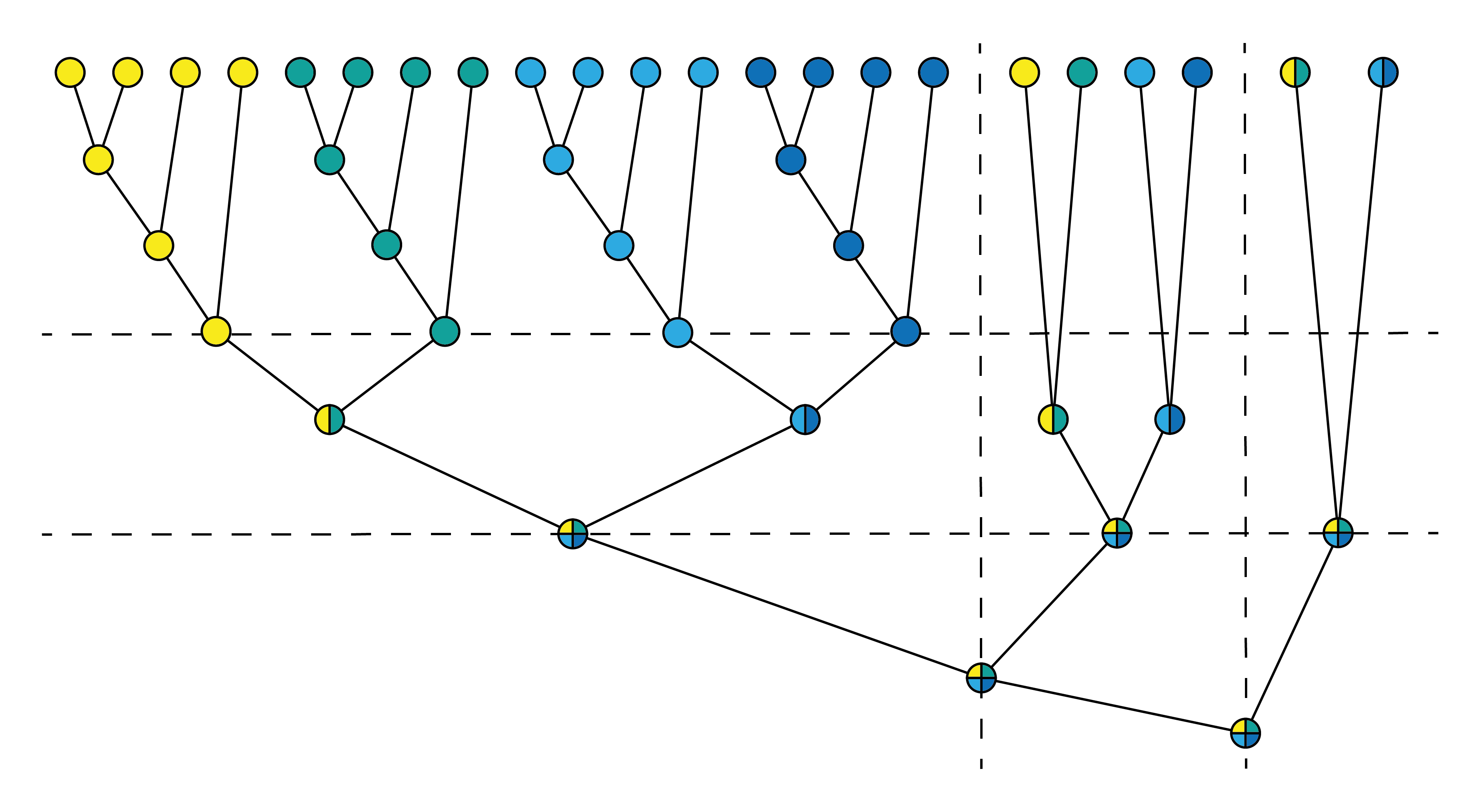}
    \caption{Illustration of the update technique as binary-tree using
      \cref{alg:welfords-weighted-online-update} on four parallel
      processes each represented with another color and over three
      estimation rounds separated by the vertical dashed lines.}
    \label{fig:triple-layered-welford}
\end{figure}

In particular as
illustrated in \cref{fig:triple-layered-welford}, \Cref{alg:welfords-weighted-online-update} is
first used to incrementally update the sample statistics on individual
processes, then again to merge the computations recursively across
multiple processing units and lastly,
\cref{alg:welfords-weighted-online-update} is utilized one more time
to update the statistical quantities over several estimation
rounds. This last step is denoted in~\cref{alg:bmlmc}, however, the
other two updates happen within \texttt{MS-FEM} for which we refer to
the upcoming~\Cref{sec:multi-sample-finite-element-method}.
Essentially, we combine in \texttt{MS-FEM} a finite element
parallelization with a sample distribution.
The resulting inherent parallelization of the algorithm
can be classified according to the
criteria defined in~\cite{badia2021massively, drzisga2017scheduling}
and~\cite{blazewicz2019handbook, dempster2018high, hager2010introduction}
as a \textit{dynamic and heterogeneous sample and solver parallelization in a single program multiple data framework}.
A detailed discussion is given in~\cite[Section 3.5.4]{baumgarten2023fully}.

Lastly, we present a proof on \Cref{proposition:bmlmc-proposition}.
The idea is to combine Gustafson's law~\cite{gustafson1988reevaluating}
with~\cref{theorem:mlmc-epsilon-theorem}.
Gustafson's law describes the theoretical slowdown of an
already parallelized task, if it is executed on a serial machine.
The motivation behind this law is to describe how more processing units
can be utilized to solve larger problems in the
same amount of time, i.e., to describe how well the parallelization scales weakly.
Translated to the knapsack problem,
larger means that we have used more samples and more levels in the final computation.
Thereby, we can achieve a smaller estimated RMSE with the same budget in time.
Hence, we can measure the weak scaling of the developed parallelization
by the development of the estimated error as more processing units are added.

\begin{proof}
    [Proof of \cref{proposition:bmlmc-proposition}]
    The estimates in the edge cases $\lambda_\rp = 1$ of an optimal parallelism
    $\epsilon \lesssim (\abs{\cP} \cdot \rT_{\mathrm{\rB}})^{-\delta} = \rC_{\epsilon}^{-\delta}$
    and $\lambda_\rp = 0$ of a serial execution
    $\epsilon \lesssim \rT_{\mathrm{\rB}}^{-\delta} =  \rT_{\epsilon}^{-\delta}$
    simply follow by inversion of \cref{theorem:mlmc-epsilon-theorem} for every feasible execution.

    To examine $\lambda_\rp \in (0, 1)$, we split the $\epsilon$-time
    into two parts of the program, a serial part $\lambda_{\rs}=1 - \lambda_{\rp}$
    and a parallel part $\lambda_{\rp}$
    executed on $\abs{\cP}$ processing units, i.e.,
    \begin{align*}
      \rT_{\epsilon, \rp} = \lambda_{\rs} \rT_{\epsilon, \rp} + \lambda_{\rp}\rT_{\epsilon, \rp}
        \,.
    \end{align*}
    If the same program is executed on a serial system,
    the parallelizable part of the system slows down by a factor of $\abs{\cP}^{-\delta}$ to
    achieve the same RMSE tolerance of $\epsilon$, i.e., the corresponding
    sequential execution time is 
    \begin{align*}
      \rT_{\epsilon, \rs} =
      \lambda_{\rs} \rT_{\epsilon, \rp} + \lambda_{\rp}\rT_{\epsilon, \rp}\abs{\cP}^{-\delta}
      \,.
    \end{align*}
    By this we can deduce the optimal speedup factor of the parallelization
    \begin{align*}
      \rS \coloneqq \frac{\rT_{\epsilon, \rs}}{\rT_{\epsilon, \rp}}
      = \frac{\lambda_{\rs} + \lambda_{\rp} \abs{\cP}^{-\delta}}{\lambda_{\rs} + \lambda_{\rp}}
      = (1 - \lambda_{\rp}) + \lambda_{\rp} \abs{\cP}^{-\delta}.
    \end{align*}
    Using this speedup factor to determine the additional error reduction by
    utilizing $\abs{\cP}$ processing units, we get again by the inverted estimate of~\cref{theorem:mlmc-epsilon-theorem}
    \begin{align*}
        \epsilon &\lesssim \rS \cdot \rT_{\mathrm{\rB}}^{-\delta}
        = \roundlr{(1 - \lambda_{\rp})  + \lambda_{\rp} \abs{\cP}^{- \delta}} \rT_{\mathrm{\rB}}^{-\delta}.
    \end{align*}
    As a consequence, there is a part in the error, denoted with $\epsilon_{\rs}$, which can only be
    reduced with further processing time and another part, denoted with $\epsilon_{\rp}$,
    which can also be mitigated by more processing units
    \begin{align*}
      \epsilon \lesssim \underbrace{(1 - \lambda_{\rp})
        \cdot \rT_{\mathrm{\rB}}^{-\delta}}_{\eqqcolon \epsilon_{\rs}}
      + \underbrace{\lambda_{\rp} \big(\abs{\cP} \cdot \rT_{\mathrm{\rB}}\big)^{-\delta}}_{\eqqcolon \epsilon_{\rp}}
      \,.
    \end{align*}
\end{proof}

We conclude this section, by summarizing the following limits as a consequence of \Cref{proposition:bmlmc-proposition}.
\begin{table}[H]
    \centering
    \label{tab:example}
    \begin{tabular}{|c|ccc|}
        \hline
        & $\lambda_{\rp} = 0$      & $\lambda_{\rp} \in (0, 1)$            & $\lambda_{\rp} = 1$      \\
        \hline
        $\rT_{\rB} \rightarrow \infty$ & $\epsilon \rightarrow 0$ & $\epsilon \rightarrow 0$              & $\epsilon \rightarrow 0$       \\
        $\abs{\cP} \rightarrow \infty$ & $\epsilon \sim 1$        & $\epsilon \rightarrow \epsilon_{\rs}$     & $\epsilon \rightarrow 0$              \\
        \hline
    \end{tabular}
\end{table}
By this tabel, the BMLMC method is
MSE-consistent with respect to the time-budget
$\rT_{\mathrm{\rB}}$ for a fixed set of processing units $\cP$. For a fixed time-budget
$\rT_{\mathrm{\rB}}$, the BMLMC method is not MSE-consistent with
respect to the amount of processing units~$\abs{\cP}$.
Hence, there remains a parallelization bias $\epsilon_{\rs}$ no matter
how manny processing units are added.
We refer to~\Cref{sec:mlmc-acoustic-experiments} for numerical experiments
on the parallelization and the derived bound.

    \section{Multi-Sample Finite Element Method}\label{sec:multi-sample-finite-element-method}

A Finite Element Method~(FEM) searches an approximation $\bu_\ell$ to
some PDE in a finite dimensional function space~$V_\ell$. To construct
this space and implement FEMs on a parallel computer,
the spatial domain $\cD \subset \RR^{D}$ is partitioned into subdomains $\cD_P$
each assigned to a different processing unit $P \in \cP$ and decomposed
in finitely many cells $K \in \cK_P$, i.e.,
\begin{equation}
    \label{eq:domain-decomposition}
    \overline{\cD}
    = \bigcup_{P \in \cP}
    \overline{\cD}_P = \bigcup_{P \in \cP} \bigcup_{K \in \cK_P}
    \overline{K}
    \quad \text{with} \quad K \cap K' = \emptyset \quad \text{for} \quad K \neq K'
    \,,
\end{equation}
where $K \subset \cD$ are open sets, $\cK = \bigcup_{P \in \cP} \cK_P$
is the collection of all cells on all processing units and
$\cK_P$ is the collection of cells on a single processing unit $P$.
The cardinality of the set of processing units $\cP$ is assumed to be
of power of two to keep the theory aligned with our implementation;
however, this is not a necessity and sets of another size might be considered as well.
This decomposition of the domain $\cD$
further defines a set of vertices $\cV$, a set of faces $\cF$ and a set of edges $\cE$
as explained in the following.

We denote with $\cF_K$ the set of faces for a cell $K \in \cK$ and for
all inner faces $F\in \cF_K\cap \cD$,
$K_F \in \cK$ represents the neighboring cell such that $\overline F =
\partial K \cap \partial K_F$.  We denote the unit normal vector on
the face $F\in \cF_K$ pointing outwards of $K$ by~$\bn_{K}$.

Furthermore, $\cV_K$ denotes the vertices of the cell $K$ and
$\cE_F$ denotes the edges of a face $F \in \cF$.
Hence, we set $\cV = \bigcup_{K \in \cK} \cV_K$, $\cF = \bigcup_{K \in \cK} \cF_K$
and $\cE = \bigcup_{F\in   \cF} \cE_F$ and define a distributed finite element mesh as
\begin{align*}
  \cM_\cP \coloneqq \set{\cV, \cK, \cF, \cE} \quad \text{with} \quad \cM_P \coloneqq \set{\cV_P, \cK_P, \cF_P, \cE_P}
  \,,
\end{align*}
where $\cV_P = \bigcup_{K\in   \cK_P} \cV_K$, $\cF_P = \bigcup_{K\in \cK_P} \cF_K$
and $\cE_P = \bigcup_{F\in   \cF_P} \cE_F$
are the vertices, faces and edges on a single processing unit.
Note that $\big\{\cV_P\big\}$, $\big\{\cF_P\big\}$, $\big\{\cE_P\big\}$ are overlapping, and
$\big\{\cK_P\big\}$ is non-overlapping for conforming discretizations; otherwise, the overlap
depends on the finite element method. 

To construct a mesh hierarchy, the cell diameter $h_0$
of a given mesh $\cM_{\ell=0, \cP}$ is sequentially
divided in half $h_\ell = h_0 2^{-\ell}$
with $\ell = 0, \dots, L$ as discretization level. This gives the hierarchy
\begin{equation*}
  \cM_{\ell=0, \cP} \subset
  \cM_{\ell=1, \cP} \subset \dots \subset \cM_{\ell=L, \cP}
  \,.
\end{equation*}
More details and several applications of this parallel data structure are given in~\cite{baumgarten2021parallel}.

\subsection{Multi-Mesh Parallelization}\label{subsec:multi-mesh-parallelization}

To combine the parallelization technique discussed
in \Cref{subsec:parallelization-techniques} with the FE parallelization,
we proceed as follows. We distribute the computational units $\cP$
across both the set of input samples $\{\by^{(m)}\}_{m=1}^{M_\ell}$ and the domain $\cD$.
Formally, this is expressed as the following resource allocation problem.

\begin{problem}
    Approximate $M_\ell$-times a PDE with a FEM on the discretization level $\ell$,
    such that the communication on a fixed set of processing units
    $\cP$ is minimized.
\end{problem}

\begin{figure}
    \centering
    \includegraphics[width=0.8\textwidth]{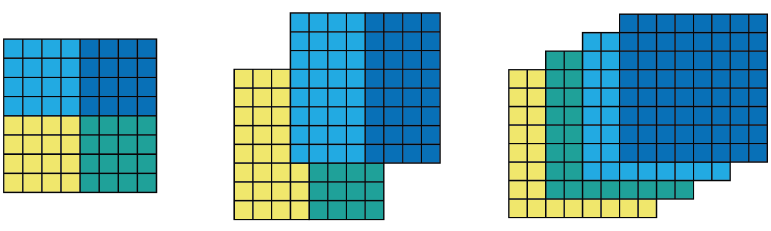}
    \caption{Parallelization approaches for $\cD = (0,1)^2$ and
        $\abs{\cP} = 4$ according
        to~\eqref{eq:comm-split-formula}. First, only on spatial domain
        $\cD$ for $M_{\ell} = 1$ with $k = 2$. Second, mixed parallelization for
        $M_{\ell} = 2$ with $k = 1$. Last, only for $M_{\ell} = 4$ with $k = 0$.}
    \label{fig:parallelization}
\end{figure}

We have to solve this problem within \Cref{alg:bmlmc} whenever the routine \texttt{MS-FEM} is invoked.
Considering that we have varying sample sizes across different levels,
we find the solution to this problem by examining the following cases.

First, we
consider the task to approximate a single $M_{\ell} = 1$ sample solution on multiple
processing units $\abs{\cP} > 1$. Then, the best parallelization is given by
the domain decomposition~\eqref{eq:domain-decomposition} resulting in a single, parallelized
mesh $\cM_\cP$ over the domain $\cD$. An illustration of this
case is given in~\cref{fig:parallelization} on the very left for
$\abs{\cP} = 4$. Second, if the sample amount equals the amount of
processing units $M_{\ell} = \abs{\cP}$, a minimal communication, and thus an
optimal parallelization, is achieved by assigning each process its
very own unparallelized mesh. This results in the set of meshes
$\cM_{\cP} \coloneqq \{\cM^{(m)}_P\}_{m=1}^{M_\ell}$, i.e.,
an individual mesh for every single sample as shown in
\cref{fig:parallelization} on the very right.

The more general case,
where we compute more samples than available processing units $M_\ell > \abs{\cP}$,
requires a sequential split of the samples $\sum_{\texttt{j}} {M}_{\texttt{j}, \ell}'
= M_\ell$ with \hbox{${M}_{\texttt{j}, \ell}' \leq \abs{\cP}$.}
Last, we consider the case as
depicted in the middle of \cref{fig:parallelization}, where $1 < M_\ell <
\abs{\cP}$. Here, we construct for each $m = 1, \dots, M_\ell$ a subset of
processing units $\cP_k^{(m)} \in \cP$ which can be used to distribute the
domain $\cD$ on. The subsets $\cP_k^{(m)}$ are disjoint and of size
$|\cP_k^{(m)}| = 2^k$ where $k \in \NN_0$ is chosen such that
\begin{equation}
    2^k \leq \frac{\abs{\cP}}{M_\ell} < 2^{k+1}.
    \label{eq:comm-split-formula}
\end{equation}
By following this rule we construct the set
of meshes $\cM_{\cP} \coloneqq \{\cM_{\cP_k}^{(m)}\}_{m=1}^{M_\ell}$,
such that we minimize the communication in every estimation round.

\subsection{Multi-Sample Finite Element Spaces}\label{subsec:multi-sample-finite-finite-spaces}

We consider the task to compute the FE solution $\bu_\ell(\omega)
\in V_{\ell}$ of $\cL_{\ell}(\omega, \bx) \bu_{\ell}(\omega, \bx) =
\bb_{\ell}(\omega, \bx)$ for multiple samples at once. In
particular, the parallel data structure presented in the previous
\Cref{subsec:multi-mesh-parallelization}
is exploited to define a finite element space incorporating
the subsets $\cP_k^{(m)} \subset \cP$.

\begin{definition}
    \label{def:multi-sample-finite-element-space}
    We call the space
    \begin{align*}
        V_{\ell}(\cP) = V_{\ell}(\cP_k^{(1)}) \times \dots \times V_{\ell}(\cP_k^{(M_\ell)})
        = \prod_{m=1}^{M_\ell} V_{\ell}(\cP_k^{(m)})
    \end{align*}
    a \textit{multi-sample finite element space}, where
    \begin{align*}
        V_\ell(\cP_k^{(m)}) \coloneqq \set{\bv_\ell \in V_{\ell} \colon \,
        \bv_\ell |_K \in V_{\ell,K}, \, \forall K \in \cK(\cP_k^{(m)}), \, \cP_k^{(m)} \subset \cP}
    \end{align*}
    is a finite element space for a single sample, defined on the triangulation $\cK(\cP_k^{(m)})$
    where the subdomain of processes $\cP_k^{(m)} \subset \cP$ is chosen with the
    rule~\eqref{eq:comm-split-formula} and $V_{\ell,K}$ is a generic local finite element space.
\end{definition}

With this definition the task of the multi-sample finite element method (MS-FEM)
is to find the coefficients
\begin{align*}
    \bmu = (\bmu_1^{(1)}, \dots, \bmu_{N_{\ell}^{h}}^{(1)}, \dots, \bmu_1^{(M_\ell)}, \dots,
    \bmu_{N_{\ell}^{h}}^{(M_\ell)})^\top \in \RR^{M_\ell \cdot N_{\ell}^{h}}
\end{align*}
representing the discrete solution
\begin{align*}
    (\bu_\ell)_{m=1}^{M_\ell} = \roundlr{\sum_{n=1}^{N_{\ell}^{h}} \bmu_n^{(m)} \bpsi_n^{(m)} }_{m=1}^{M_\ell}  \in V_{\ell}(\cP) \,,
\end{align*}
where $\bpsi_n^{(m)}$ are basis functions of the global finite element space
of dimension $N_{\ell}^{h}$.

This formulation is inspired by the implementation
in~\cite{wieners2023mpp315}, where the parallelization over the
samples is realized on the coefficient vector of the FEM.
This enables the highly adaptive parallelization scheme
needed in the BMLMC method. Finally, the complete procedure
is summarized in~\cref{alg:ms-fem}.

\begin{algorithm}
    \caption{Multi-Sample Finite Element Method}
    \label{alg:ms-fem}
    \vspace{-2mm}
    \begin{align*}
        &\texttt{function MS-FEM}(M_{\ell}, \cP, \ell) \colon \\
        &\begin{cases}
             k \leftarrow \texttt{use } \eqref{eq:comm-split-formula} \texttt{ with } (M_\ell, \cP) \\
             \texttt{Generate input samples } \set{\by^{(m)}}_{m=1}^{M_\ell} \\
             \texttt{Assemble } \cL_{\ell}(\by^{(m)}, \bx) \, \bu_\ell(\by^{(m)}, \bx) = \bb_\ell(\by^{(m)}, \bx)
             \texttt{ for } m = 1, \dots, M_\ell \\
             \texttt{Solve } \cL_{\ell}(\by^{(m)}, \bx) \, \bu_\ell(\by^{(m)}, \bx) = \bb_\ell(\by^{(m)}, \bx)
             \texttt{ for } m = 1, \dots, M_\ell \\
             \texttt{Update QoI and Cost with \cref{alg:welfords-weighted-online-update}}
        \end{cases}
    \end{align*}
\end{algorithm}

\begin{remark}
    The above is applicable to arbitrary finite element spaces,
    e.g.~continuous Lagrange elements, enriched Galerkin elements,
    Raviart-Thomas elements, space-time discontinuous Galerkin (dG) elements
    or, as in the upcoming section, to dG elements in space. Further details
    and experiments can be found in~\cite{baumgarten2023fully}.

    Since the load distribution is a function of $M_\ell$ and $\abs{\cP}$
    minimizing the communication, the
    system in~\Cref{alg:ms-fem} is assembled, such that it minimizes the coupling, i.e.,
    the system is decoupled for each sample and mildly coupled on the
    spatial domain. The assembled system has a block structure and is sparse which is inherited from the
    sparsity of each finite element discretization block.
\end{remark}

    \section{Discretization of the Acoustic Wave Equation}
\label{sec:discretization-of-the-acoustic-wave-equation}

In our numerical examples, we consider the acoustic wave
equation with randomly modeled input data in
the form of compressible waves propagating through solids.

\begin{problem}
    \label{problem:stochastic-acoustic}
    Let $\cD \subset \RR^D$ be a domain and
    $[0, T] \subset \RR$ a time interval.
    We search for the randomly distributed velocity field $\bv \colon \Omega \times \cD \times [0, T] \rightarrow \RR^{D}$
    and pressure component $p \colon \Omega \times \cD \times [0, T] \rightarrow \RR$, such that
    \begin{align*}
        \pdeProblem{
            \rho(\omega, \bx) \partial_{t} \bv(\omega, \bx, t) - \nabla p(\omega, \bx, t)
            &=& \mathbf{f}(\omega, \bx, t) &\bx \in \cD, \,\,\,\,\, t \in (0, T] \\
            \kappa(\omega, \bx)^{-1} \partial_{t} p(\omega, \bx, t) - \div \roundlr{\bv(\omega, \bx, t)}
            &=& g(\omega, \bx, t) &\bx \in \cD, \,\,\,\,\, t \in (0, T] \\
            \bv(\omega, \bx, t) \cdot \bn
            &=& 0 &\bx \in \partial \cD, \,\, t \in [0, T] \\
            \mathbf{v}(\omega, \bx, 0) &=& \mathbf{v}_0(\omega, \bx) &\bx \in \cD \\
            p(\omega, \bx, 0) &=& p_0(\omega, \bx)  &\bx \in \cD
        }
    \end{align*}
    with $\mathbf{f} \colon \Omega \times \cD \times (0, T]
    \rightarrow \RR^{D}$ and $g \colon \Omega \times \cD \times
    (0, T] \rightarrow \RR$ as right-hand sides and the material
    parameters $\kappa, \rho \colon \Omega \times \cD
    \rightarrow \RR$ modeled as random fields. We further allow
    for randomly distributed initial data in the velocity
    component $\bv_0 \colon \Omega \times \cD \rightarrow
    \RR^{D}$ and the pressure component $p_0 \colon \Omega
    \times \cD \rightarrow \RR$.
\end{problem}

The works~\cite{mishra2012sparse, mishra2016multi} discuss sufficient conditions for the
well-posedness of~\cref{problem:stochastic-acoustic} in the case $\cD = \RR^D$.

\subsection{Semi-Discretization with Discontinuous Galerkin Methods}
\label{subsec:semi-discretization-with-discontinuous-galerkin-methods}

We follow~\cite{bohlen2020visco} and use a discontinuous Galerkin approximation in space
based on the formulation of the acoustic wave equation as a first-order system for a fixed $\omega \in \Omega$, $\bx \in \cD$
and $t \in (0, 1]$ given by
\begin{equation}
    \label{eq:first-order-formulation}
    \rM(\omega, \bx) \partial_t \bu(\omega, \bx, t) + \rA \bu(\omega, \bx, t)
    = \bb(\omega, \bx, t) \quad \text{and} \quad
    \bu(\omega, \bx, 0) = \bu_0(\omega, \bx)
    \,.
\end{equation}
From now on, we omit the explicit notation of the dependency on $\omega$, $\bx$ and $t$.
The first oder formulation~\eqref{eq:first-order-formulation} is derived with the operators
\begin{align*}
    \bu =
    \begin{pmatrix}
        \bv \\  p
    \end{pmatrix},
    \quad
    \rM \bu =
    \begin{pmatrix}
        \rho \, \bv \\ \kappa^{-1} p
    \end{pmatrix}, \quad
    \rA \bu = -
    \begin{pmatrix}
        \nabla p \\ \div \bv
    \end{pmatrix}
    \quad \text{and} \quad
    \bb =     \begin{pmatrix}
                  \mathbf{f} \\ g
    \end{pmatrix}.
\end{align*}
This system is approximated in space using discontinuous Galerkin (dG) finite elements
\begin{align*}
    V_{\ell, \bp}^\mathrm{dG} =
    \set{\bu_\ell \in \rL^{2}(\cD; \RR^{D+1}) \colon \bu_\ell |_K \in V_{K, \bp}^\mathrm{dG}, \, \forall K \in \cK}
    \,,
\end{align*}
where $V_{K, \bp}^\mathrm{dG} \coloneqq \QQ_\bp(K; \RR^{D+1})$ is the
tensor product space of local polynomials on a cell $K \in \cK$.
In the resulting semi-discrete system, we search for $\bu_\ell \in V_{\ell, \bp}^\mathrm{dG}$
\begin{equation*}
    \rM_\ell \partial_t \bu_\ell + \rA_\ell \bu_\ell = \bb_\ell
    \quad \text{and} \quad \bu_\ell(0) = \bu_{\ell,0}
    \,,
\end{equation*}
with cell-wise constant approximations for $\rM_\ell$ and $\rL_2$-projections of $\bb, \bu_{0}$ on
$\bb_\ell,\bu_{\ell,0} \in V_{\ell, \bp}^\mathrm{dG}$.
The differential operator $\rA$ is discretized with a full-upwind scheme
\begin{align*}
    \sprod{\rA_{\ell} \bu_\ell, \bphi_\ell}_{\rL^2(\cD)}
    = \sum_{K\in \cK} \sprod{\rA_{\ell, K} \bu_\ell, \bphi_{\ell,K}}_{\rL^2(K)} \,,
\end{align*}
with test functions $\bphi_\ell = (\bvarphi_\ell,\psi_\ell) \in V_{\ell, \bp}^\mathrm{dG}$.
Each local operator is given in case of Neumann boundary conditions by
\begin{align*}
    \sprod{\rA_{\ell, K} \bu_\ell, \bphi_{\ell,K}}_{\rL^2(K)} =
    &- \sprod{\nabla p_{\ell, K}, \bvarphi_{\ell,K}}_{\rL^2(K)}
    - \sprod{\div \bv_{\ell, K}, \psi_{\ell,K}}_{\rL^2(K)}
    \\
    &\hspace{-2.7cm}- \sum_{F \in \cF_K \cap \cD} \frac{1}{Z_K + Z_{K_F}}
    \sprod{\bbracketlr{p_{\ell, K}}_F + Z_{K_F} \bbracketlr{\bv_{\ell, K}}_F \cdot \bn_{K}, \, \psi_{\ell, K}
        + Z_K \bvarphi_{\ell, K} \cdot \bn_{K}}_{\rL^2(F)}
    \\
    &\hspace{-2.7cm} + \sum_{F \in \cF_K\cap\partial\cD} Z_K
    \sprod{\bv_{\ell, K} \cdot \bn_K, \, \psi_{\ell, K}
        + Z_K \bvarphi_{\ell, K} \cdot \bn_K}_{\rL^2(F)}
\end{align*}
where $Z_K = \sqrt{\kappa_K \rho_K}$ is the impedance and
$\bbracketlr{\bv_{\ell, K}}_F = \bv_{\ell}|_{K_F} - \bv_{\ell}|_{K}$
is the jump at inner faces $F \in \cF \cap \cD$. In~\cite{corallo2022space}
it is shown that the system is the well-posed,
also for more general boundary conditions.

\subsection{Time-Discretization with Implicit Methods}
\label{subsec:time-discretization-with-implicit-methods}

We follow~\cite[Section 3]{bohlen2020visco} and shortly outline the
usage of the implicit mid-point rule with the time-step size
$\tau_{\ell} = T / N_{\ell}^{\tau}$ and the time-steps $t_n = n \tau_{\ell}$, $n=0,
\dots, N_{\ell}^{\tau}$, i.e., we construct a sequence of approximations
$\bu_\ell(t_n) \in V_{\ell, \bp}^{\text{dG}}$ with the initial value
of $\bu_\ell(t_0) = \bu_0$ by
\begin{equation}
    \label{eq:acoustic-implicit-mpr}
    \roundlr{\rM_\ell + \frac{\tau_{\ell}}{2} \rA_\ell} \bu_\ell(t_n)
    = \roundlr{\rM_\ell - \frac{\tau_{\ell}}{2} \rA_\ell} \bu_\ell(t_{n - 1}) + \tau_{\ell} \bb_{\ell}(t_{n - 1/2}).
\end{equation}

By~\cite[Theorem 3.1]{bohlen2020visco} the above system is well-posed
and therefore the implicit midpoint rule is applicable. The downside
of this implicit method is that the new iteration $\bu^{n+1}$ is only
implicitly given. Thus, a system of algebraic equations has
to be solved in each time-step which can increase the cost of the
method significantly.

However, in the context of hyperbolic PDEs the
usage of implicit methods avoids stability issues, if the
Courant–Friedrichs–Lewy (CFL) condition~\cite{courant1928partiellen}
\begin{equation}
    \label{eq:cfl}
    \frac{\tau_{\ell}}{h_\ell} \leq C_{\mathrm{CFL}} \quad \Longleftrightarrow \quad \tau_{\ell} \leq C_{\mathrm{CFL}} \cdot h_0 2^{-\ell}
\end{equation}
is not satisfied. For our particular problem, this is critical since
the wave speed is a random variable in each cell and thus, the
right ratio $\tau_{0} / h_0 = C_{\mathrm{CFL}}$ is a local condition
leading to global stability issues. The usage of implicit methods
avoids this problem, nevertheless, finding the right ratio is still
important since too large time-steps lead to worse conditioned systems
in~\eqref{eq:acoustic-implicit-mpr} and too many time-steps simply
might be unnecessary to achieve a smaller overall error.
However, it is shown in~\cite[Lemma~3.1]{bohlen2020visco}
that~\eqref{eq:acoustic-implicit-mpr} is well-conditioned
and the convergence is independent of the mesh size on level $\ell$.
In particular, we solve this system using a GMRES solver
with a point block Jacobi preconditioner.
We refer to \Cref{sec:mlmc-acoustic-experiments} for an experimental
investigation of this issue.

    \section{Numerical Experiments for the Acoustic Wave Equation}
\label{sec:mlmc-acoustic-experiments}

We consider~\cref{problem:stochastic-acoustic} and solve it with the
methods introduced in the previous sections. In particular, we commit
to the following problem and method configurations which will serve,
if not stated otherwise, as the default for the numerical
experiments.

\begin{figure}[H]
    \centering
    \begin{adjustbox}{raise=0.55cm}
        \includegraphics[width=0.17\textwidth]{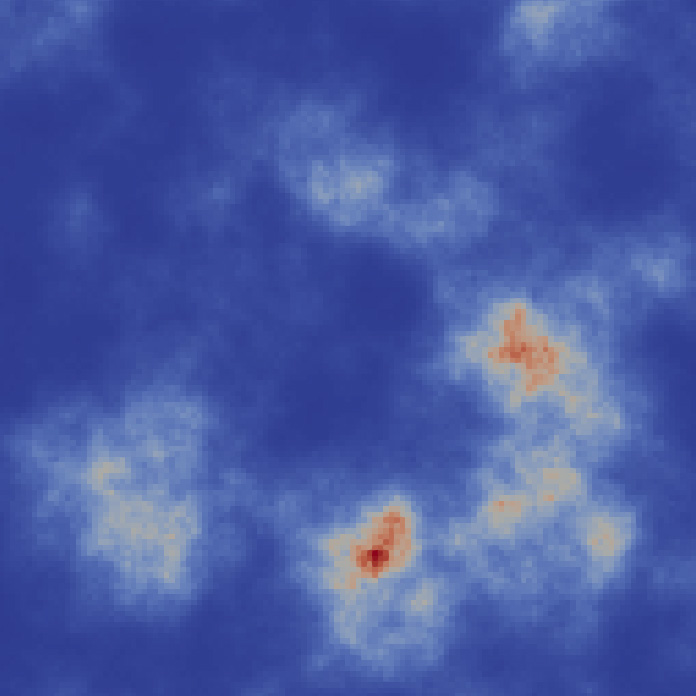}
    \end{adjustbox}
    \hspace{0.1cm}
    \includegraphics[width=0.25\textwidth]{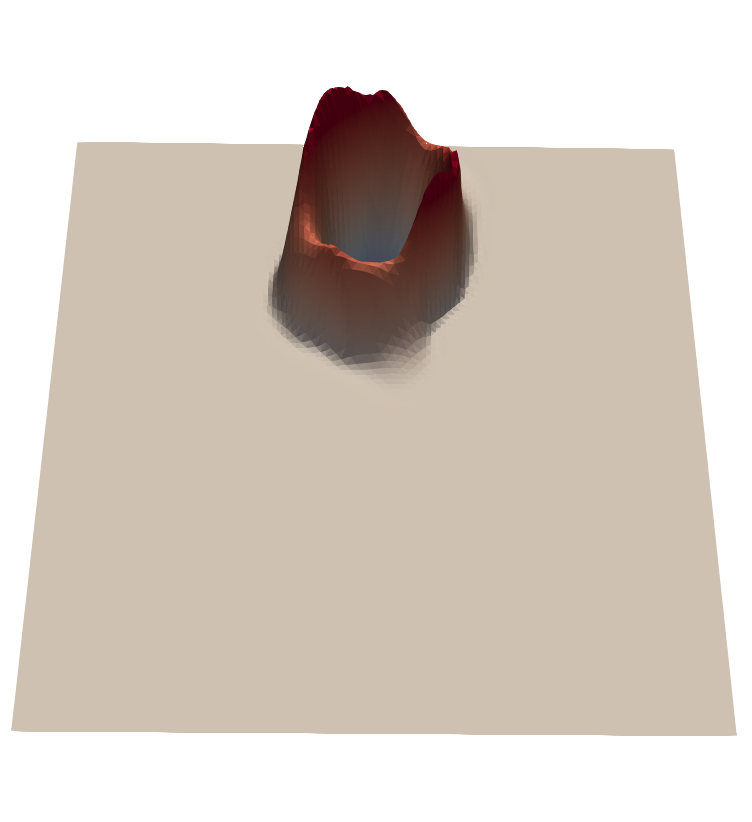}
    \includegraphics[width=0.25\textwidth]{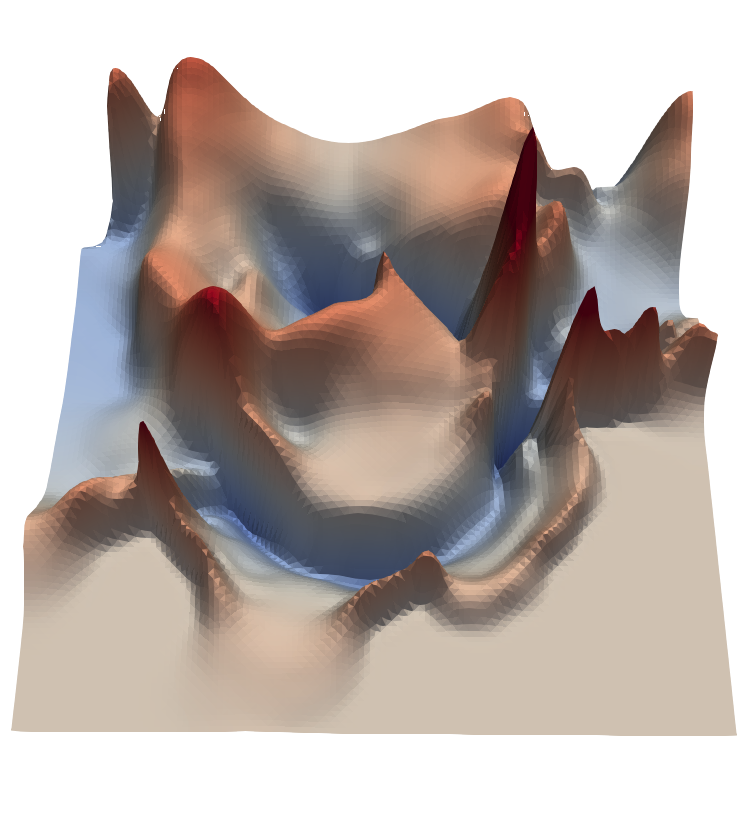}
    \includegraphics[width=0.25\textwidth]{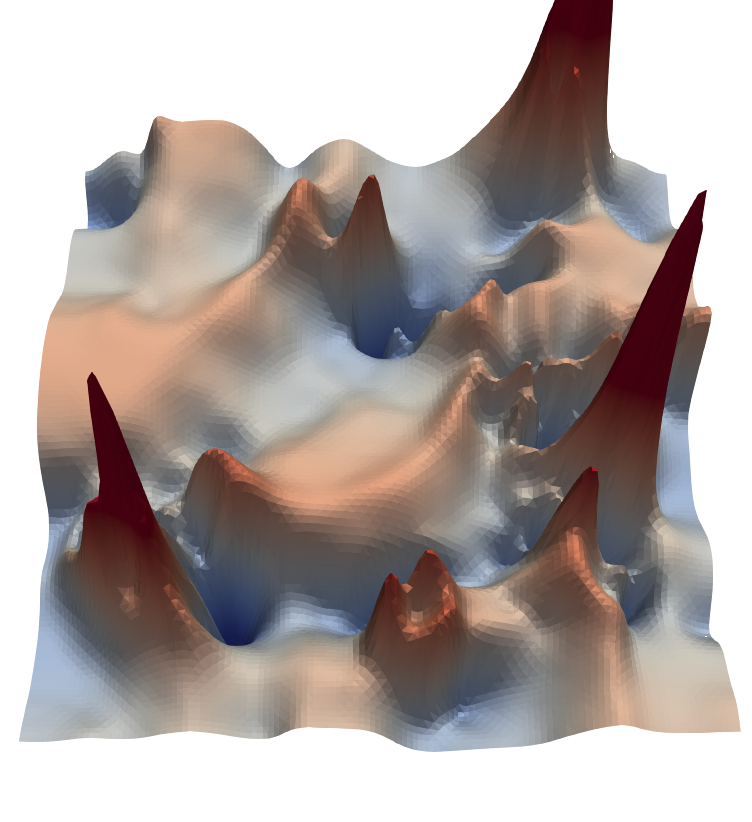}
    \vspace{-5mm}
    \caption{One example realization of the pressure wave
    approximating~\cref{problem:stochastic-acoustic} with log-normal
    material density $\rho$ (left image) shown at the time points
        $t=0.0625$, $t=0.5$ and $t=1.0$. The regions of high material
        density (red areas in left plot) lead to small wave speeds in
        these regions and thus, to high pressures in the wave.}
    \label{fig:time-stepping-plots}
\end{figure}

\paragraph{Problem Configuration}

We consider the domain $\cD = (0, 1)^2$, the final time $T = 1$ and
homogeneous and deterministic initial conditions $(\bv_0, p_0) =
(\mathbf{0}, 0)^\top$.

The right-hand side of \cref{problem:stochastic-acoustic} is deterministically
given by $\mathbf{f} \equiv \mathbf{0}$ and $g(\bx, t) = g_1(t) \,
g_2(\bx)$, where the function $g_1(t)$ is a Ricker wavelet, i.e.,
\begin{align*}
    g_1(t) &= 10 \roundlr{1 - \roundlr{\tfrac{t}{a}}^2} \cdot \exp \roundlr{- \tfrac{t^2}{2 a^2}}
    \quad \text{with} \quad a = \tfrac{\pi}{10}, \quad t \in [0, 1].
\end{align*}
The function $g_2(\bx)$ is a nascent delta function centered at $\bc =
(0.5, 0.75)^{\top}$ with an appropriate constant $\overline{g}_2$ such
that $\norm{g_2}_{\rL^1(\cD)} = 1$
and a diameter $w = 0.1$, i.e.,
\begin{align*}
    g_2(\bx) &=
    \begin{cases}
        \overline{g}_2 \exp \roundlr{- \roundlr{1 - \norm{\frac{\bx - \bc}{w}}_2^2}^{-1}}, &\norm{\bx - \bc}_2 < w \\
        0, &\norm{\bx - \bc}_2 \geq w
    \end{cases} \qquad \bx \in \cD.
\end{align*}
As material, we use a uniformly constant and deterministic compression module $\kappa \equiv 1$ and
a log-normally distributed material density $\rho(\omega, \bx)$, i.e., $\log(\rho(\omega, \bx))$ is
Gaussian random field with mean-zero and the covariance function
\begin{equation}
    \label{eq:covariance-function}
    \Cov(\bx_1, \bx_2) = \sigma^2 \exp \roundlr{- \norm{\frac{\bx_1 - \bx_2}{\lambda}}_2^{\nu}},
\end{equation}
where the variance $\sigma = 1.0$, the correlation length $\lambda = 0.15$
and the smoothing $\nu = 1.8$ are used. Defining
\begin{equation*}
    \rho_{\min}(\omega) \coloneqq \min_{\bx \in \cD} \rho(\omega, \bx)
    \quad \text{and} \quad \rho_{\max}(\omega) \coloneqq \max_{\bx \in \cD} \rho(\omega, \bx)
\end{equation*}
gives a distribution of the maximal
$c_{\max}(\omega) = \sqrt{\kappa / \rho_{\min}} < \infty$
and the minimal $c_{\min}(\omega) = \sqrt{\kappa / \rho_{\max}} > 0$ wave speeds.
By~\cite[Lemma 2.3]{charrier2013finite} realizations are Hölder continuous,
and thus  $\rho(\omega, \cdot) \in \rL^\infty(\cD)$ for a fixed $\omega \in \Omega$.
Lastly, we mention that the samples of $\rho$ are
generated with the circulant embedding method~\cite{dietrich1997fast}
on the multi-mesh implementation introduced
in~\ref{subsec:multi-mesh-parallelization}. For further details
we refer to~\cite{baumgarten2023fully}.

The default QoI is the $\rL^2$-norm for vector valued functions
in a region of interest $\cD^{\text{RoI}} = (0.25, 0.75) \times (0, 0.25)$ at time $T = 1$, i.e.,
\begin{align*}
    \rQ(\omega) \coloneqq \roundlr{\int_{0.25}^{0.75} \int_{0}^{0.25}
        \abs{(\bv, p)^\top(\omega, x_1, x_2, 1)}_2^2 \rd x_2 \, \rd x_1}^{1/2},
\end{align*}
where $\abs{\cdot}_2$ is the Euclidean-norm.
The problem configuration is illustrated in~\cref{fig:time-stepping-plots}
for one particular realization of the input data at the time points $t=0.0625$, $t=0.5$ and $t=1.0$.

\paragraph{Method Configuration}
The experiments are conducted on the HoReKa
supercomputer for $\rT_{\rB} = 6$ hours using $\abs{\cP} = 1024$ processing units.
We initialize the BMLMC method on four initial levels starting with the mesh width $h_0 = 2^{-5}$ by
\begin{align*}
    \set{M_{0, \ell}^{\text{init}}}_{\ell=0}^{L_0} \coloneqq \set{M_{0, 0}^{\text{init}} = 2^{12},
        \, M_{0, 1}^{\text{init}} =  2^{10}, \, M_{0, 2}^{\text{init}} = 2^{7}, \, M_{0, 3}^{\text{init}} = 2^{5}}
\end{align*}
which consumes less than 5\% of the total computational budget but already provides good initial estimates.
Furthermore, we choose the splitting factor as $\theta = 0.5$ and the reduction factor as $\eta = 0.9$.
The semi-discrete solution is
searched in $V_{\ell, \bp=2}^{\text{dG}}$ on uniform meshes with $h_\ell = h_0 2^{-\ell}$,
which is then solved using
an implicit midpoint rule with the time-step size $\tau_\ell = 2^{-\ell + 3}$.

\paragraph{Covariance Function}

\begin{figure}
    \centering
    \includegraphics[width=1.0\textwidth]{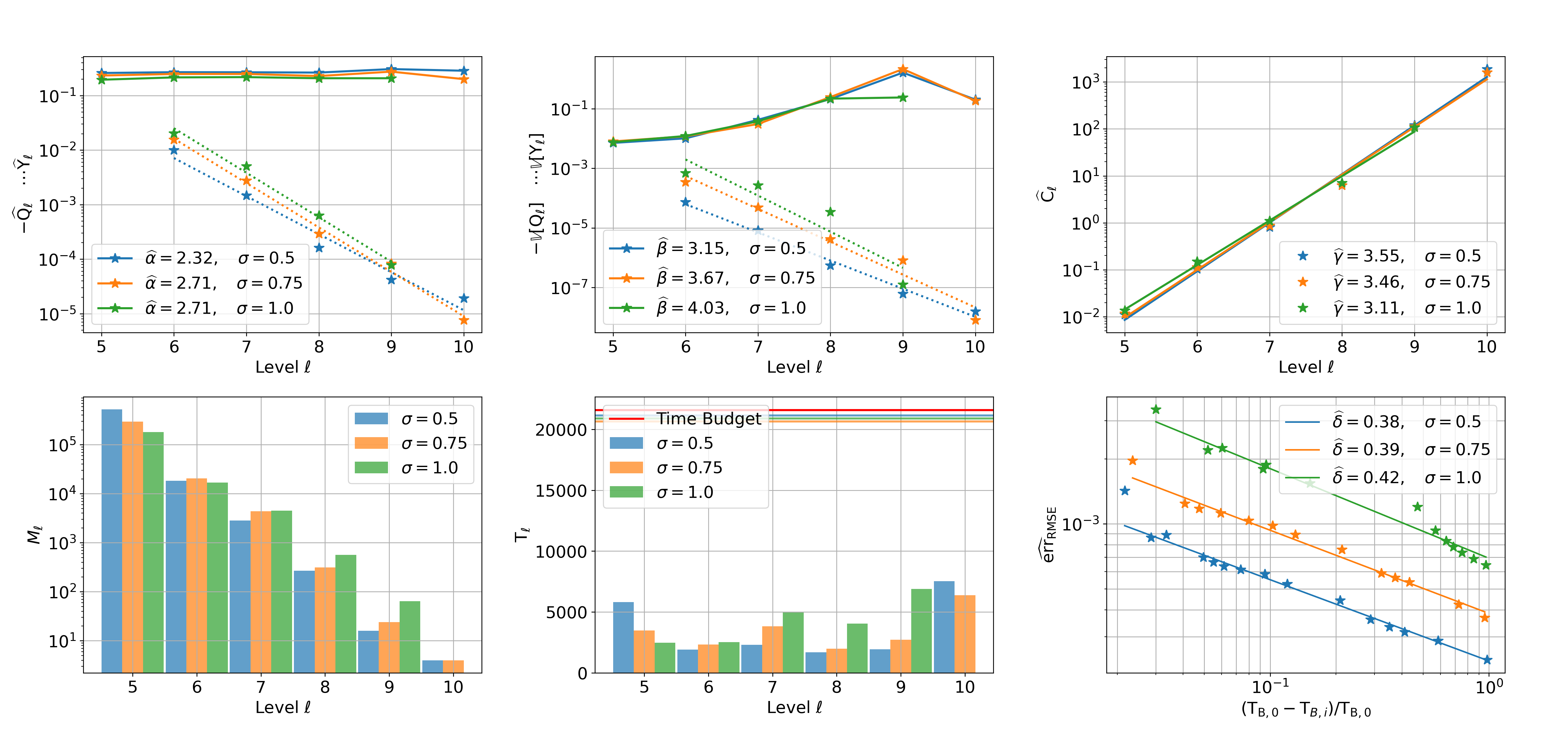}
    \vspace{-10mm}
    \caption{Experiments on the covariance function~\eqref{eq:covariance-function} of the log-normal field.}
    \label{fig:mlmc-acoustic-time-stepping-sigma-on-horeka}
\end{figure}

As start, we examine the influence of the covariance
function~\eqref{eq:covariance-function} on the behavior
of~\cref{alg:bmlmc}. Analytical investigations
~\cite{charrier2012strong, charrier2013finite, nobile2015multi, teckentrup2013further}
as well as experiments~\cite{baumgarten2023fully, baumgarten2021parallel} for
elliptic problems have shown that the structure of the log-normal
fields has a large influence on the constant and the convergence rate
in~\cref{theorem:mlmc-epsilon-theorem}. We conduct similar
investigations for the acoustic wave equation with log-normally
distributed material parameters by choosing $\sigma \in \set{0.5,
    0.75, 1.0}$ in the covariance
function~\eqref{eq:covariance-function}, while everything else is kept
as described in the con\-figurations. The results of this experiment are
given in~\cref{fig:mlmc-acoustic-time-stepping-sigma-on-horeka}, where
in the top row the a posteriori verification of~\cref{assumption:mlmc}
is given with the estimated exponents $\widehat{\alpha},
\widehat{\beta}$ and $\widehat{\gamma}$. In the bottom row, the figure
shows the computed amount of samples on each level
on the left, the cost distribution over the
levels in the middle, and the numerical verification of the
convergence of~\cref{proposition:bmlmc-proposition} on the right. The x-axis
of the lower right plot is the relative left over time budget
$(\rT_{\rB, 0} - \rT_{\rB, \ttti})/\rT_{\rB, 0}$ and the y-axis is the
estimated RMSE over the estimation rounds in logarithmic
scales. \cref{fig:mlmc-acoustic-time-stepping-sigma-on-horeka} clearly
shows that increasing the variance in~\eqref{eq:covariance-function}
worsens the constant in~\cref{proposition:bmlmc-proposition}, while the
measured convergence rate $\widehat{\delta}$, estimated by
\begin{align*}
    \min_{(\widehat{\delta}, \, \widehat{c}_{\delta})}
    \sum_{\ttti} \big( \log_2 (\widehat{\err}_{\text{RMSE}, \ttti}) - \widehat{\delta}
    \log_2((\rT_{\rB, 0} - \rT_{\rB, \ttti}) / \rT_{\rB, 0}) + \widehat{c}_{\delta} \big)^2
    \,,
\end{align*}
only changes slightly. We
further remark that the BMLMC method works very reliably for this
model problem and is capable to exhaust the large computational budget
of $\rB = 1024 \cdot 6$ CPU hours feasibly and completely. This can be
seen on the middle plot on the bottom, where the total computing times
are given by the horizontal lines staying just below the
time represented by the red line. Similar
investigations for $\nu$ and $\lambda$
in~\eqref{eq:covariance-function} or any sort of input data
to~\cref{problem:stochastic-acoustic} can be done as well for which we
refer again to~\cite{baumgarten2023fully}.

\paragraph{Time Discretization}

\begin{figure}
    \centering
    \includegraphics[width=1.0\textwidth]{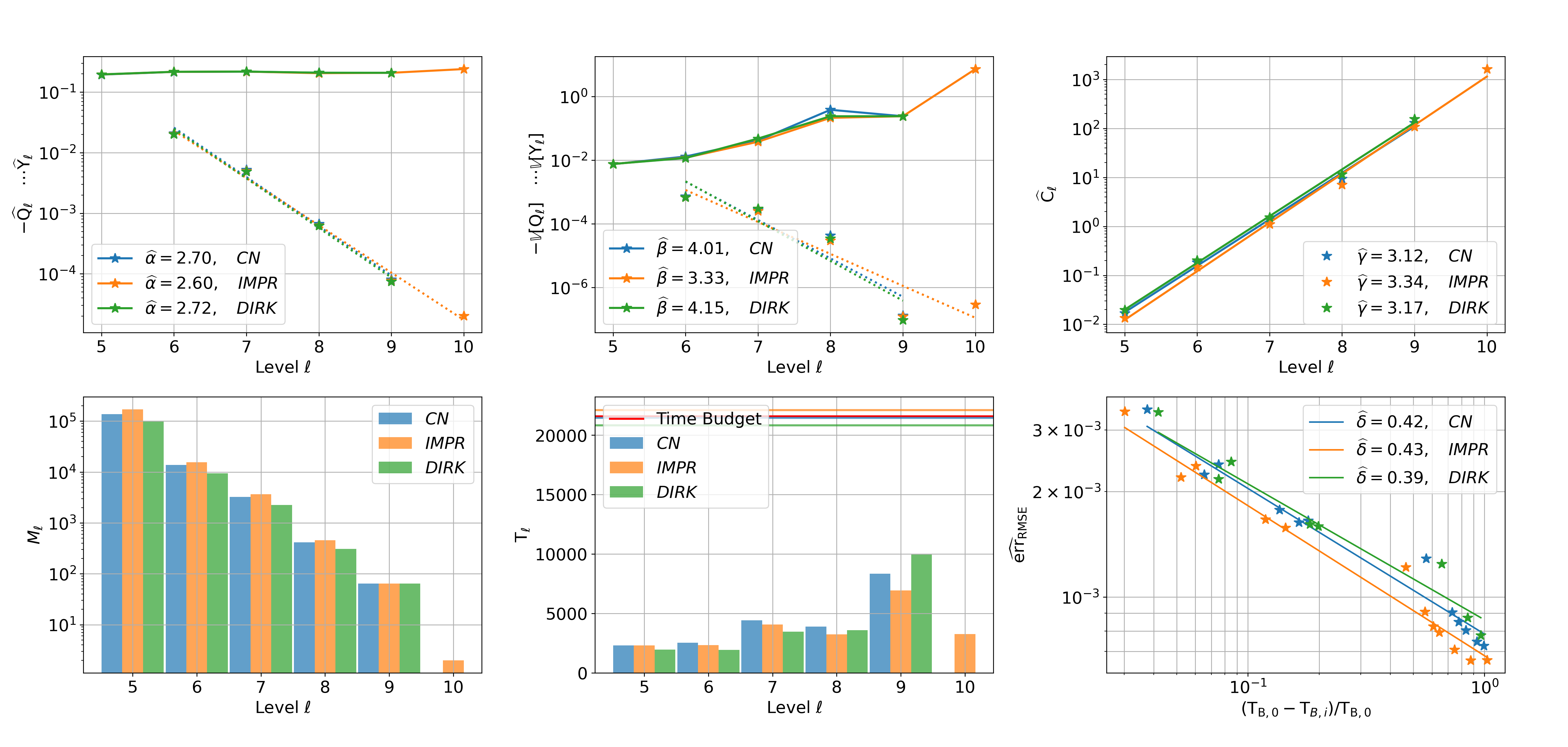}
    \vspace{-10mm}
    \caption{Comparison of different time-stepping methods.}
    \label{fig:mlmc-acoustic-time-stepping-rkorder-on-horeka}
\end{figure}

We further investigate the time discretization. Even though only
shortly discussed in~\Cref{subsec:time-discretization-with-implicit-methods},
finding the right time-steps and the right time integrator is crucial
for the performance of the overall method and its
stability. In~\cref{fig:mlmc-acoustic-time-stepping-rkorder-on-horeka},
we illustrate the comparison of three different implicit Runge-Kutta
methods with the global convergence order of
$\cO(\tau_{\ell}^2)$. Particularly, we compare the implicit midpoint rule
(IMPR), the Crank Nicolson (CN) method and a third diagonal implicit
Runge-Kutta (DIRK) method determined by the Butcher-tableau:
\begin{align*}
    \text{DIRK} \colon \, \begin{array}
    {c|cc}
                              1/4 & 1/4 & 0   \\
                              3/4 & 1/2 & 1/4 \\
                              \hline
                              & 1/2 & 1/2
    \end{array}
\end{align*}

By the lower right plot of
\cref{fig:mlmc-acoustic-time-stepping-rkorder-on-horeka}, we see that
the implicit midpoint rule yields the smallest estimated error and
thereby is the best choice out of theses three since we have assigned
all three experiments the same computational budget. We suspect that
the reason for this is that the evaluations in each time-step in the
IMPR are cheaper than for the other two methods. As a consequence of
this cost saving, more samples and even one additional level can be
computed using the IMPR. We remark that we experimented with explicit
Runge-Kutta methods, too, but the time-step sizes $\tau_\ell$ had to
be drastically reduced in order to stabilize the computation. Locally
adaptive schemes as in~\cite{grote2021uncertainty} might overcome this
issue, however, we have not been comparing this ansatz to the current
implicit approach yet.

\begin{figure}[H]
    \centering
    \includegraphics[width=1.0\textwidth]{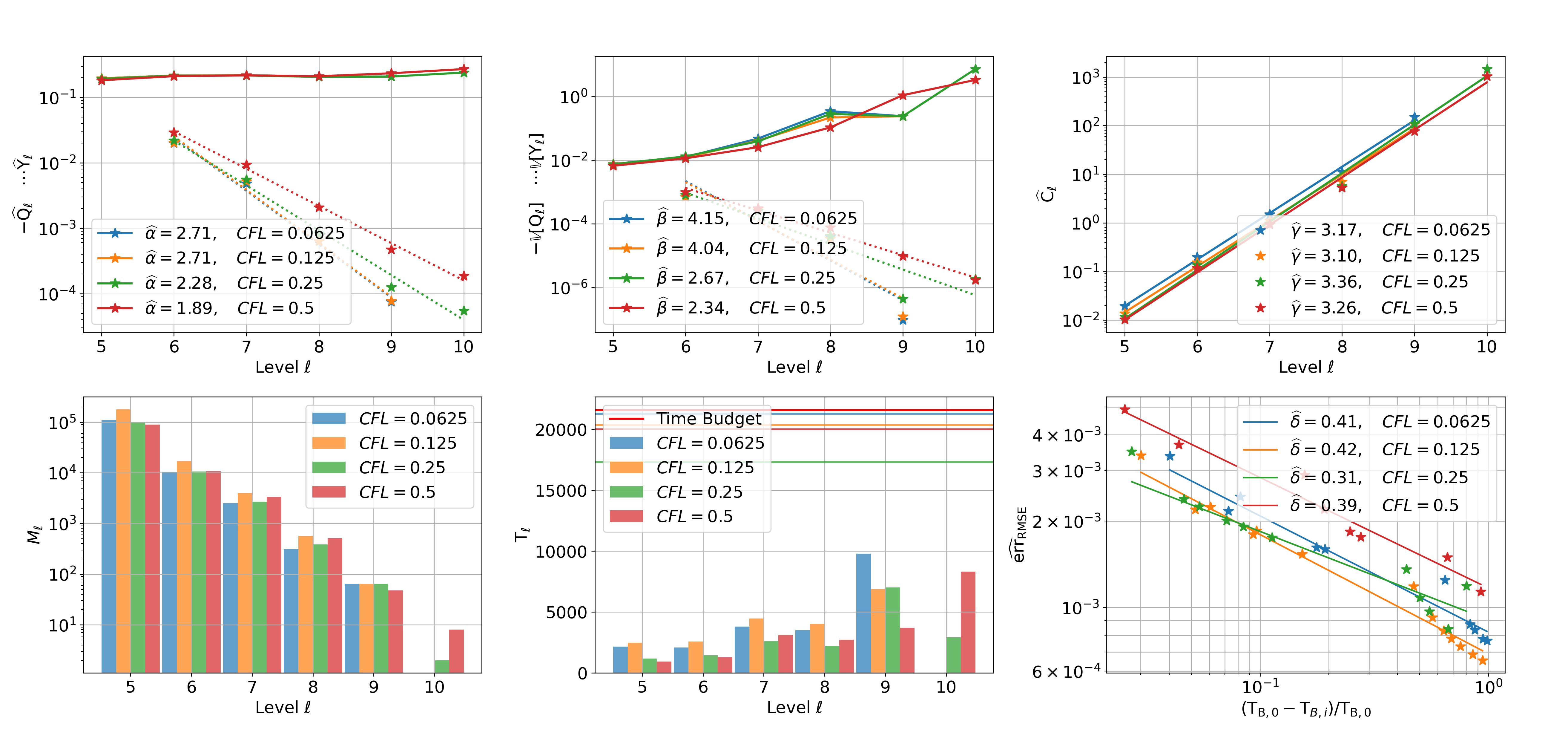}
    \vspace{-10mm}
    \caption{Comparison of different time-step to mesh-width ratios $C_{\mathrm{CFL}} = \tau_\ell / h_\ell$.}
    \label{fig:mlmc-acoustic-time-stepping-cfl-on-horeka}
\end{figure}

We recall the discussion of~\Cref{subsec:time-discretization-with-implicit-methods} and examine the
influence of the time-step size on the overall method performance. The
results are given
in~\cref{fig:mlmc-acoustic-time-stepping-cfl-on-horeka} where we tried
out different ratios $\tau_{\ell} / h_{\ell} \in \set{2^{-1}, 2^{-2},
    2^{-3}, 2^{-4}}$. By the plot on the upper right, we see that the
constant $c_{\gamma}$ slightly depends on the time-step as predicted,
but also that the variance reduction $\widehat{\beta}$ (upper row in
the middle) is heavily influenced. The best choice is
$C_{\mathrm{CFL}} = 2^{-3}$ again reviled in the lower right plot
of~\cref{fig:mlmc-acoustic-time-stepping-cfl-on-horeka}. We further
remark that with this choice the estimate $\widehat{\gamma} = 3.1$ is
very close to the best possible value of $\gamma = D + 1 = 3$ as
further explained in~\cite{hochbruck2015efficient}.

\paragraph{Space Discretization}

\begin{figure}
    \centering
    \includegraphics[width=1.0\textwidth]{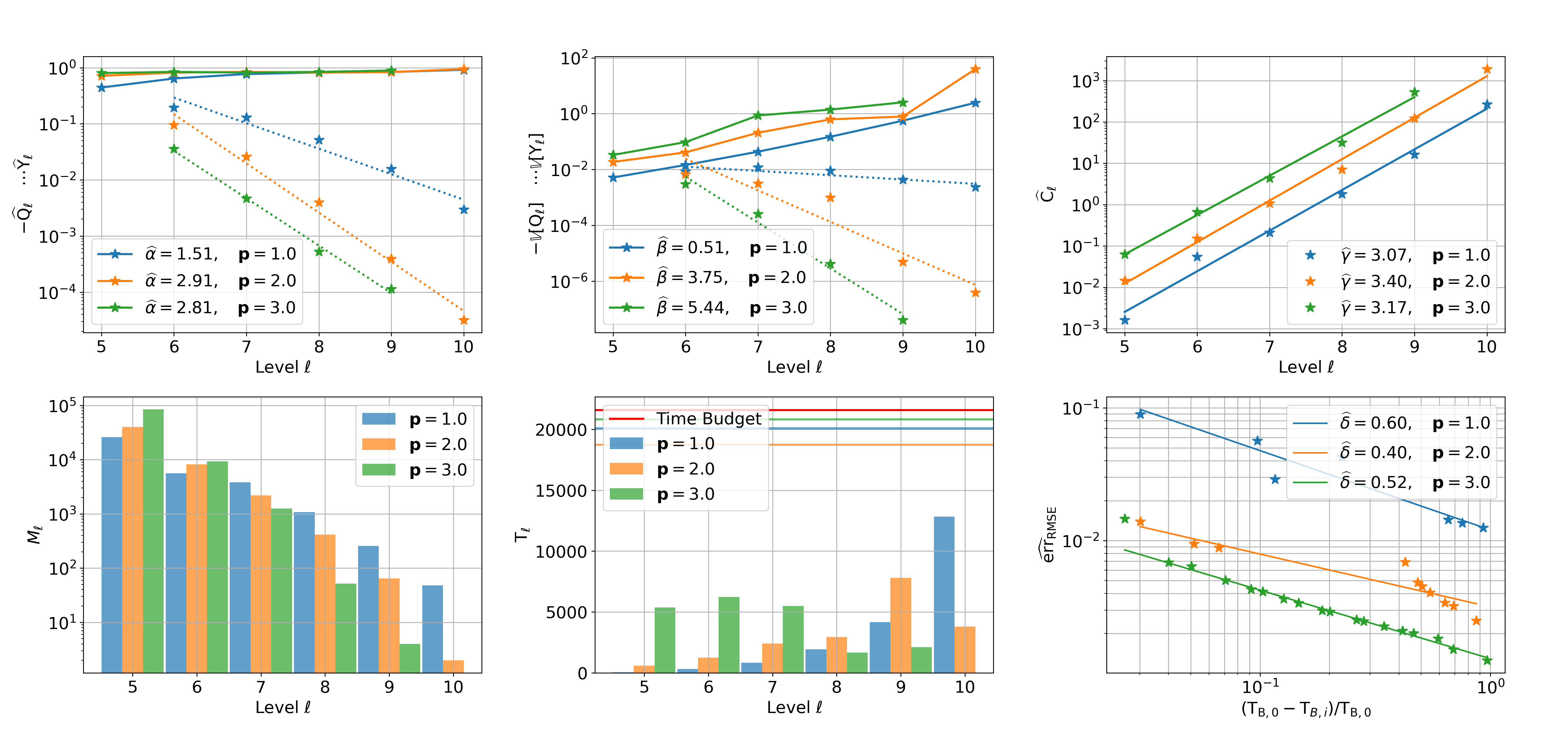}
    \vspace{-10mm}
    \caption{Comparison of different polynomial degrees in the ansatz space $V_{\ell, \bp}^{\rd\rG}$.}
    \label{fig:mlmc-acoustic-time-stepping-degree-on-horeka}
\end{figure}

For the next experiment, we are interested in the polynomial degree of
the dG space $V_{\ell, \bp}^{\rd\rG}$. The results
in~\cref{fig:mlmc-acoustic-time-stepping-degree-on-horeka} show that
an ansatz space with a higher degree is worth to consider since the
$\widehat{\err}_{\mathrm{RMSE}}$ gets smaller with a growing degree
even though the cost constant $c_{\gamma}$ (confer upper right plot)
is higher. We emphasize that this conclusion is highly problem
dependent and that the higher polynomial degree is only worth the
additional cost, if the true solution to the PDE provides enough
regularity. It is well known, for
example given in a discussion in~\cite{cliffe2011multilevel}, that the
cost is dominated by the highest level if $\beta > \gamma$. Contrary
to that, if $\beta < \gamma$, the cost is dominated by the lower
levels. Both cases can be observed
in~\cref{fig:mlmc-acoustic-time-stepping-degree-on-horeka} on the bar
plot in the center of the bottom row, where for $\bp = 1$ the
estimated exponents satisfy $\widehat{\beta} = 0.51 < 3.07 =
\widehat{\gamma}$ and for $\bp = 3$ the exponents are given as
$\widehat{\beta} = 5.44 > 3.17 = \widehat{\gamma}$.

\paragraph{Parallelization}

\begin{figure}
    \centering
    \includegraphics[width=1.0\textwidth]{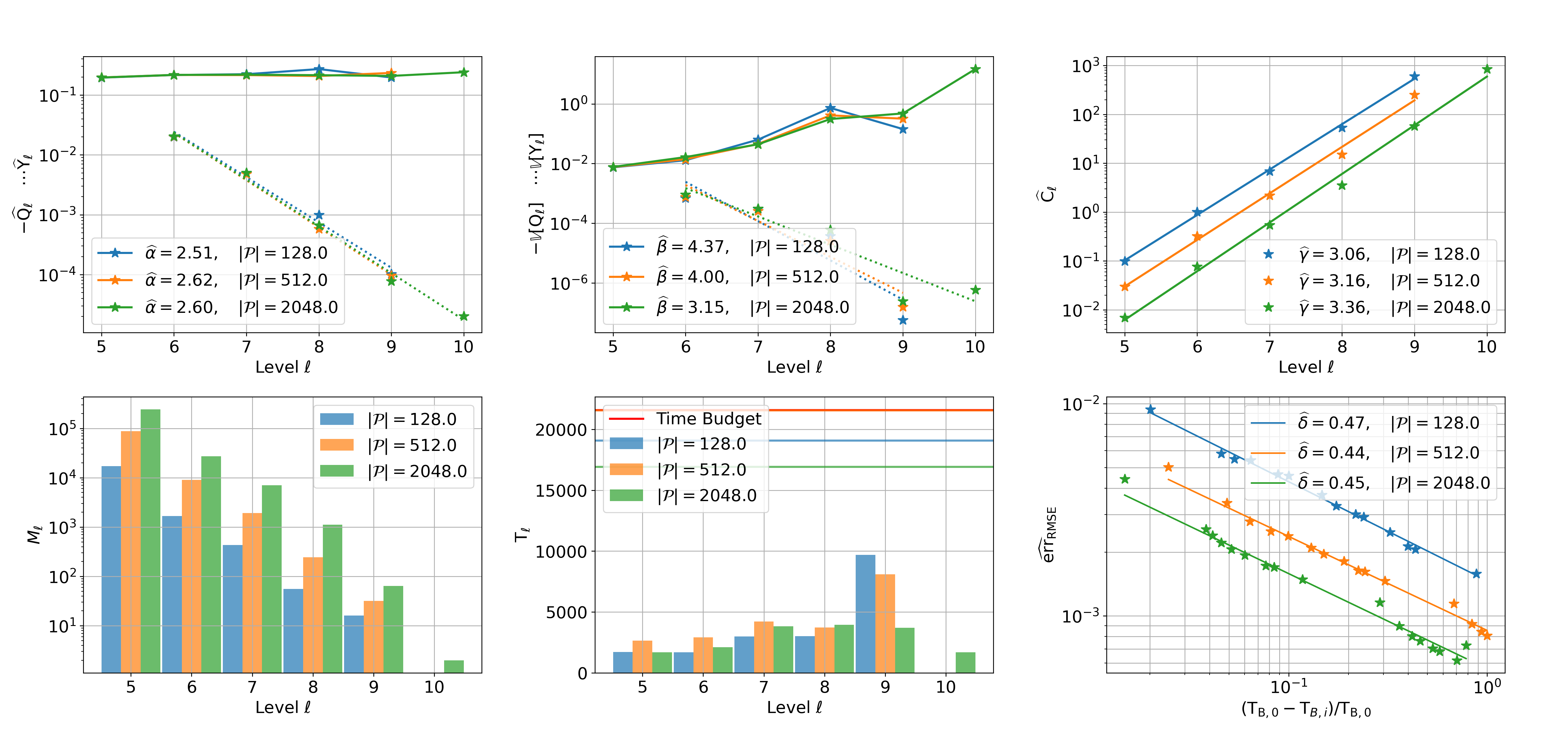}
    \vspace{-10mm}
    \caption{Weak scaling parallelization experiment for a fixed time budget $\rT_{\rB}$.}
    \label{fig:mlmc-acoustic-time-stepping-weak-scaling-on-horeka-2x3}
\end{figure}

Last but not least, we examine the proposed parallelization by conducting a weak
scaling experiment, i.e., we increase the computational resources from
$\abs{\cP}=128$ to $\abs{\cP}=2048$ and keep the computational time
budget fixed at $\rT_{\rB} = 6$ hours. The numerical results of this experiment are
summarized in~\Cref{fig:mlmc-acoustic-time-stepping-weak-scaling-on-horeka-2x3}.
The lower right plot indicates clearly that we effectively reduce the estimated
error by utilizing more processing units. However, to examine the influence
of the reduction factor $\eta$ and to evaluate the method
in the light of \Cref{proposition:bmlmc-proposition},
we solve the problem again on $\abs{\cP_{\max}}=8192$ and on $\abs{\cP_{\max}} \cdot 2^{-k}$
with $k=1,\dots,7$. Subsequently, we consider the estimated error
at the very end of the simulation, i.e.~at
$(\rT_{\rB, 0} - \rT_{\rB, \ttti}) / \rT_{\rB, 0} = 1$
and plot this over $2^{-k}$.
With this and~\eqref{eq:bmlmc-theorem}, we conclude for some $\lambda_{\rp} \in (0, 1)$
\begin{align*}
    \epsilon_{k} \lesssim \epsilon_{\rs} + \epsilon_{\rp, k}
    = \epsilon_{\rs} + \lambda_{\rp} (2^{-k} \cdot \abs{\cP_{\max}} \cdot \rT_{{\rB}})^{-\delta}
    = \epsilon_{\rs} + \lambda_{\rp} (\abs{\cP_{\max}} \cdot \rT_{{\rB}})^{-\delta} \cdot 2^{k \delta}.
\end{align*}

\begin{figure}[h]
    \centering
    \includegraphics[width=0.5\textwidth]{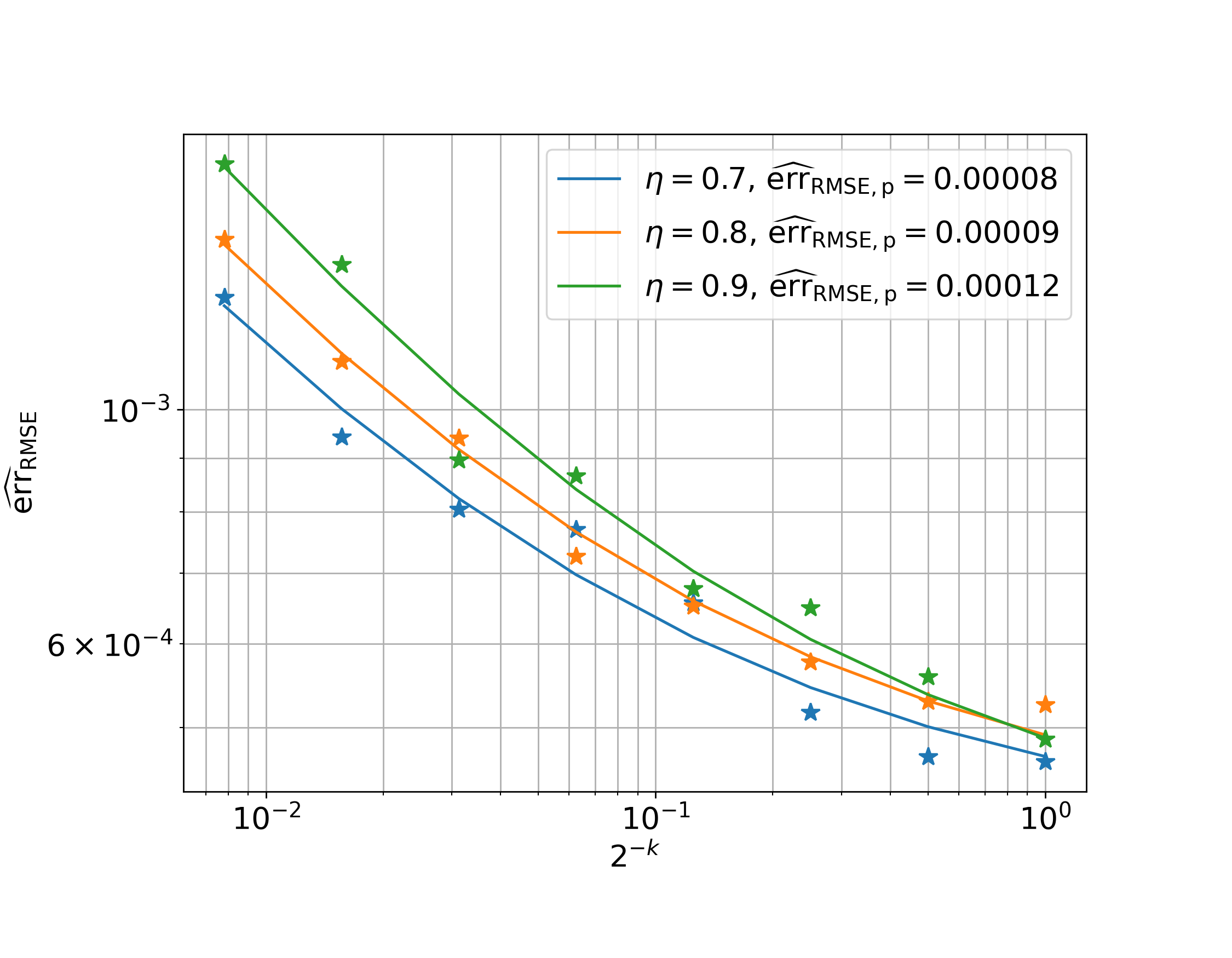}
    \vspace{-7mm}
    \caption{Numerical verification of \Cref{proposition:bmlmc-proposition}.}
    \label{fig:mlmc-acoustic-time-stepping-weak-scaling-on-horeka}
\end{figure}

This motivates to determine $\widehat{\err}_{\text{RMSE}, \rs}$
and $\widehat{\err}_{\text{RMSE}, \rp}$ by fitting the curve
\begin{align*}
    \widehat{\err}_{\text{RMSE}, k} = \widehat{\err}_{\text{RMSE}, \rs}
    + \widehat{\err}_{\text{RMSE}, \rp} \cdot 2^{k \widehat{\delta}}.
\end{align*}
The results of these experiments and the fitted curve are given in
\Cref{fig:mlmc-acoustic-time-stepping-weak-scaling-on-horeka} for $\eta \in \set{0.7, 0.8, 0.9}$,
which illustrates the influence of the reduction factor on the parallelization.
Clearly, we can see in this plot that the smaller the reduction factor,
the smaller is the estimated error. This is because the larger reduction factor
leads to more frequent synchronizations of the processing units and thus,
leading to parallelization losses and ultimately in larger errors.
The downside of small reduction factors is the higher probability
of exceeding or not fully using the computational budget.
In conclusion, in \Cref{fig:mlmc-acoustic-time-stepping-weak-scaling-on-horeka}
we see that the proposed BMLMC method adheres to the theoretical bound of
\Cref{proposition:bmlmc-proposition} and that $\lambda_{\rp}$ is mostly
influenced by the reduction factor $\eta$.

    \section{Discussion, Conclusion and Outlook}\label{sec:discussion-and-outlook}

We present a novel adaptation of the MLMC method called Budgeted MLMC (BMLMC) method.
This approach minimizes the need for prior knowledge, while demonstrating high reliability,
robustness, and broad applicability. Furthermore, it achieves exceptional
performance within budget constraints and offers full parallelization
up to the limits of Gustafson's law.

The method's effectiveness stems from three fundamental components:
the seamless integration of MLMC with FE methods, the adaptive load distribution
within a single distributed data structure according to~\eqref{eq:comm-split-formula},
and the resource allocation within an HPC system following
the optimality principle~\eqref{eq:bmlmc-bellman}.

To demonstrate this experimentally, we conduct investigations on the challen\-ging problem of approximating
acoustic wave equations in random and heterogeneous media. Our methodology involves a
fully automated process using the continuous delivery pipeline of the
software M++~\cite{wieners2023mpp315} connected to the HoReKa supercomputer.
This allows us to reproduce and
enhance the numerical results obtained from our implementation.

For a comprehensive explanation of the software, we refer to
a forthcoming publication or to~\cite{baumgarten2023fully, baumgarten2021parallel}.
This will provide
detailed insights into the software development workflow, as well as highlight
the distinguishing features and applications of~M++ including space-time
discretizations~\cite{corallo2022space, dorfler2016space},
interval arithmetic computations~\cite{wunderlich2022computer},
full waveform inversion~\cite{bohlen2020visco},
and other challenging applications like cardio-vascular
simulations~\cite{froehlich2022cardiac}.

This utilization of automated investigations and the empirical search for the
optimal algorithm combination provide motivation to view the MLMC method as a
knapsack problem. This perspective has facilitated the development of the BMLMC method,
incorporating DP techniques and drawing inspiration from the
continuation MLMC method~\cite{collier2015continuation}.

The same approach could potentially be applied to other multi-level UQ
algorithms such as multi-level stochastic collocation (MLSC)
or multi-level quasi-Monte Carlo (MLQMC) methods. To begin, we recognize that
we can utilize the same parallelization strategy and the same distributed data structure.
However, DP relies on a nested problem structure. By employing nested
sparse grids and suitable lattice rules for QMC, we can preserve this problem
structure. This results in a less flexible selection of
$\{{\vartriangle} M_{\ttti, \ell}\}_{\ell=0}^{L_\ttti}$.
Furthermore, the on-the-fly estimation of the errors is not as
straightforward as it is for simple Monte Carlo methods.
Although these methods may be more complex,
we anticipate that the benefits of
the BMLMC approach can be extended to MLSC and MLQMC methods with careful
consideration and adaptation. We leave detailed investigations
into these directions open for future work.

Finally, we interpret the consequence of \Cref{proposition:bmlmc-proposition}.
The notion that infinite computing power does not bypass computing time is derived by
combining the $\epsilon$-cost theorem and Gustafson's law. Intuitively, this result
makes sense because adding more workers also entails increased synchronization,
which prevents the objective from being arbitrarily optimized.
Thus, there will always be a synchronization or parallelization bias
that can only be mitigated by allocating more time for optimization.
The key aspect of demonstrating this insight lies in viewing computing time
as a component of the algorithm's cost and connecting the algorithm's objective
with hardware resources through Gustafson's law. This idea is quite general,
suggesting that similar statements to those in \Cref{proposition:bmlmc-proposition}
can be derived for other UQ algorithms beyond MLMC,
and perhaps even in the realm of machine learning.

    \section*{Acknowledgement}

We acknowledge the financial support by the CRC 1173 on Wave Phenomena,
the technical support by the National High-Performance Computing Center (NHR) at KIT
and the detailed feedback by Tobias Jahnke on an earlier draft of this work.

    \bibliographystyle{siamplain}

    \bibliography{bibliography}
\end{document}